\documentclass{emsprocart}

\contact[bartelsa@math.uni-muenster.de]{Arthur Bartels, Mathematisches
  Institut, Westf\"alische Wilhelms-Universit\"at M\"unster,
  Einsteinstr.~62, D-48149 M\"unster, Germany}

\contact[echters@math.uni-muenster.de]{Siegfried Echterhoff,
Mathematisches Institut, Westf\"alische Wilhelms-Universit\"at M\"unster,
  Einsteinstr.~62, D-48149 M\"unster, Germany}

\contact[lueck@math.uni-muenster.de]{Wolfgang L\"uck, Mathematisches
  Institut, Westf\"alische Wilhelms-Universit\"at M\"unster,
  Einsteinstr.~62, D-48149 M\"unster, Germany}

\title{Inheritance of Isomorphism Conjectures under colimits}
\author{Arthur Bartels, Siegfried Echterhoff and Wolfgang L\"uck}

\usepackage{hyperref,amscd}
\usepackage{calc}
\usepackage{enumerate,amssymb}
\usepackage[arrow,curve,matrix,tips,2cell]{xy}
  \SelectTips{eu}{10} \UseTips
  \UseAllTwocells

\DeclareMathAlphabet{\matheurm}{U}{eur}{m}{n}

\newcommand{\CastAlg}{C^*\text{-}\matheurm{Algebras}}

\newcommand{\finker}{\matheurm{finker}}

\newcommand{\Groupoids}{\matheurm{Groupoids}}

\newcommand{\Or}{\matheurm{Or}}
\newcommand{\Rings}{\matheurm{Rings}}
\newcommand{\astRings}{\ast\!\text{-}\matheurm{Rings}}
\newcommand{\Spaces}{\matheurm{Spaces}}
\newcommand{\Spectra}{\matheurm{Spectra}}

\DeclareMathOperator{\asmb}{asmb}
\DeclareMathOperator{\aut}{aut}

\DeclareMathOperator{\cent}{cent}

\DeclareMathOperator{\coker}{coker}
\DeclareMathOperator{\colim}{colim}

\DeclareMathOperator{\id}{id}

\DeclareMathOperator{\im}{im}

\DeclareMathOperator{\ind}{ind}
\DeclareMathOperator{\map}{map}
\DeclareMathOperator{\mor}{mor}

\newcommand{\pt}{\{\bullet\}}
\DeclareMathOperator{\pr}{pr}

\DeclareMathOperator{\topo}{top}

\newcommand{\Fin}{{\mathcal{F}\text{in}}}

\newcommand{\VCyc}{{\mathcal{VC}\text{yc}}}

  \newcommand{\IC}{\mathbb{C}}

  \newcommand{\IZ}{\mathbb{Z}}

  \newcommand{\calc}{\mathcal{C}}

  \newcommand{\calf}{\mathcal{F}}
  \newcommand{\calg}{\mathcal{G}}
  \newcommand{\calh}{\mathcal{H}}

  \newcommand{\calk}{\mathcal{K}}

  \newcommand{\bfE}{{\mathbf E}}
  
  \newcommand{\bff}{{\mathbf f}}
  \newcommand{\bfF}{{\mathbf F}}

  \newcommand{\bfI}{{\mathbf I}}

  \newcommand{\bfK}{{\mathbf K}}
  \newcommand{\bfKH}{{\mathbf K}{\mathbf H}}
  \newcommand{\bfL}{{\mathbf L}}

  \newcommand{\bfT}{{\mathbf T}}

\newcommand{\KH}{{K\!H}}

\newcounter{commentcounter}

\newcommand{\xycomsquare}[8]                
{\xymatrix{#1 \ar[r]^-{#2} \ar[d]^{#4} &
    #3 \ar[d]^{#5}  \\
    #6\ar[r]^-{#7} & #8 }}

\newcommand{\EGF}[2]{E_{#2}(#1)}            

\theoremstyle{plain}
\newtheorem{theorem}{Theorem}[section]

\newtheorem{lemma}[theorem]{Lemma}
\newtheorem{corollary}[theorem]{Corollary}

\theoremstyle{definition}
\newtheorem{definition}[theorem]{Definition}

\newtheorem{remark}[theorem]{Remark}

\theoremstyle{remark}

\renewcommand{\theenumi}{(\roman{enumi})}

\makeatletter\let\c@equation=\c@theorem\makeatother
\renewcommand{\theequation}{\thetheorem}

\begin{document}

\maketitle

\begin{abstract}
  We investigate when Isomorphism Conjectures, such as the ones due to
  Baum-Connes, Bost and Farrell-Jones, are stable under colimits of groups over
  directed sets (with not necessarily injective structure maps). 
  We show in particular that both the $K$-theoretic Farrell-Jones
  Conjecture and the Bost Conjecture with coefficients hold for those groups for which
  Higson, Lafforgue and Skandalis have disproved the Baum-Connes Conjecture with
  coefficients.\\[3mm]
  keywords: Bost Conjecture, inheritance under colimits.\\[2mm]
  AMS-classification: 19K99, 18F25, 55N91.
\end{abstract}


\typeout{---------------------------- Introduction ------------------------------------}
\setcounter{section}{-1}
\section{Introduction}
\label{sec:introduction}


\subsection{Assembly maps}
\label{subsec:Assembly_maps}

We want to study the following \emph{assembly maps}:

\begin{flalign}
  & \asmb^G_n \colon && H_n^G(\EGF{G}{\VCyc};\bfK_R) && \to && H_n^G(\pt;\bfK_R) =
  K_n(R\rtimes G);
\label{ass_K_n}
\\
& \asmb^G_n \colon && H_n^G(\EGF{G}{\Fin};\bfKH_R) && \to && H_n^G(\pt;\bfKH_R) =
\KH_n(R\rtimes G);
\label{ass_KH_n}
\\
& \asmb^G_n \colon && H_n^G(\EGF{G}{\VCyc};\bfL_R^{\langle -\infty\rangle}) && \to &&
H_n^G(\pt;\bfL_R^{\langle -\infty\rangle}) = L_n^{\langle - \infty \rangle}(R \rtimes G);
\label{ass_L_n}
\\
& \asmb^G_n \colon && H_n^G(\EGF{G}{\Fin};\bfK^{\topo}_{A,l^1}) && \to &&
H_n^G(\pt;\bfK^{\topo}_{A,l^1}) = K_n(A \rtimes_{l^1} G);
\label{ass_K_ntopo_l1}
\\
& \asmb^G_n \colon && H_n^G(\EGF{G}{\Fin};\bfK^{\topo}_{A,r}) && \to &&
H_n^G(\pt;\bfK^{\topo}_{A,r}) = K_n(A \rtimes_r G);
\label{ass_K_ntopo_r}
\\
& \asmb^G_n \colon && H_n^G(\EGF{G}{\Fin};\bfK^{\topo}_{A,m}) && \to &&
H_n^G(\pt;\bfK^{\topo}_{A,m}) = K_n(A \rtimes_m G).
\label{ass_K_ntopo_m}
\end{flalign}

Some explanations are in order.  A \emph{family of subgroups of $G$} is a collection of
subgroups of $G$ which is closed under conjugation and taking subgroups.  Examples are the
family $\Fin$ of finite subgroups and the family $\VCyc$ of virtually cyclic subgroups.

Let $\EGF{G}{\calf}$ be the \emph{classifying space associated to $\calf$}.  It is
uniquely characterized up to $G$-homotopy by the properties that it is a $G$-$CW$-complex
and that $\EGF{G}{\calf}^H$ is contractible if $H \in \calf$ and is empty if $H \notin
\calf$.  For more information about these spaces $\EGF{G}{\calf}$ we refer for instance to
the survey article~\cite{Lueck(2005s)}.

Given a group $G$ acting on a ring (with involution) by structure preserving maps, let $R
\rtimes G$ be the twisted group ring (with involution) and denote by $K_n(R\rtimes G)$,
$\KH_n(R\rtimes G)$ and $L_n^{\langle -\infty \rangle}(R\rtimes G)$ its \emph{algebraic
  $K$-theory} in the non-connective sense (see Gersten~\cite{Gersten(1972)} or
Pedersen-Weibel~\cite{Pedersen-Weibel(1985)}), \emph{its homotopy $K$-theory} in the sense
of Weibel~\cite{Weibel(1989)}, and its \emph{$L$-theory with decoration $-\infty$} in the
sense of Ranicki~\cite[Chapter~17]{Ranicki(1992a)}.  Given a group $G$ acting on a
$C^*$-algebra $A$ by automorphisms of $C^*$-algebras, let $A\rtimes_{l^1} G$ be the Banach
algebra obtained from $A \rtimes G$ by completion with respect to the $l^1$-norm, let $A
\rtimes_r G$ be the reduced crossed product $C^*$-algebra, and let $A\rtimes_m G$ be the
maximal crossed product $C^*$-algebra and denote by $K_n(A \rtimes_{l^1} G)$, $K_n(A
\rtimes_r G)$ and $K_n(A \rtimes_m G)$ their \emph{topological $K$-theory}.

The source and target of the assembly maps are given by $G$-homology theories (see
Definition~\ref{def:G-homology-theory} and
Theorem~\ref{the:Construction_of_equivariant_homology_theories}) with the property that
for every subgroup $H \subseteq G$ and $n \in \IZ$
$$\begin{array}{lcl} H_n^G(G/H;\bfK_R) & \cong & K_n(R\rtimes H);
  \\
  H_n^G(G/H;\bfKH_R) & \cong & \KH_n(R\rtimes H);
  \\
  H_n^G(G/H;\bfL_R^{\langle -\infty\rangle}) & \cong &
  L_n^{\langle-\infty\rangle}(R\rtimes H);
  \\
  H_n^G(G/H;\bfK^{\topo}_{A,l^1}) & \cong & K_n(A\rtimes_{l^1} H);
  \\
  H_n^G(G/H;\bfK^{\topo}_{A,r}) & \cong & K_n(A\rtimes_r H);
  \\
  H_n^G(G/H;\bfK^{\topo}_{A,r}) & \cong & K_n(A\rtimes_m H).
\end{array}$$
All the assembly maps are induced by the projection from $\EGF{G}{\Fin}$ or
$\EGF{G}{\VCyc}$ respectively to the one-point-space $\pt$.

\begin{remark} It might be surprising to the reader that we restrict to C*-algebra 
coefficients $A$ in the assembly map~\eqref{ass_K_ntopo_l1}.  Indeed, 
our main results rely heavily 
on the validity of the conjecture for hyperbolic groups, which, so far, 
is only known for C*-algebra
coefficients. Moreover we also want to study the passage from the $l^1$-setting 
to the $C^*$-setting.
Hence we decided to restrict ourselves to the case of $C^*$-coefficients throughout.
We mention that on the other hand the assembly map \eqref{ass_K_ntopo_l1} can also
be defined for Banach algebra coefficients \cite{Paravicini(2006)}.
\end{remark}


\subsection{Conventions}
\label{subsec:Conventions}

Before we go on, let us fix some conventions.  A group $G$ is always discrete.  Hyperbolic
group is to be understood in the sense of Gromov (see for instance~\cite{Bowditch(1991)},
\cite{Bridson-Haefliger(1999)}, \cite{Ghys-Harpe(1990)}, \cite{Gromov(1987)}).  Ring means
associative ring with unit and ring homomorphisms preserve units.  Homomorphisms of Banach
algebras are assumed to be norm decreasing.


\subsection{Isomorphism Conjectures}
\label{subsec:Isomorphism_Conjectures}

The \emph{Farrell-Jones Conjecture} for algebraic $K$-theory for a
group $G$ and a ring $R$ with $G$-action by ring automorphisms says
that the assembly map~\eqref{ass_K_n} is bijective for all $n \in
\IZ$.  Its version for homotopy $K$-theory says that the assembly
map~\eqref{ass_KH_n} is bijective for all $n \in \IZ$.  If $R$ is a
ring with involution and $G$ acts on $R$ by automorphism of rings with
involutions, the $L$-theoretic version of the Farrell-Jones Conjecture
predicts that the assembly map~\eqref{ass_L_n} is bijective for all $n
\in \IZ$. The Farrell-Jones Conjecture for algebraic $K$- and
$L$-theory was originally formulated in the paper by
Farrell-Jones~\cite[1.6 on page 257]{Farrell-Jones(1993a)} for the
trivial $G$-action on $R$. Its homotopy $K$-theoretic version can be
found in~\cite[Conjecture~7.3]{Bartels-Lueck(2006)}, again for trivial
$G$-action on $R$.

Let $G$ be a group acting on the $C^*$-algebra $A$ by automorphisms of $C^*$-algebras. The
\emph{Bost Conjecture with coefficients} and the \emph{Baum-Connes Conjecture with
  coefficients} respectively predict that the assembly map~\eqref{ass_K_ntopo_l1}
and~\eqref{ass_K_ntopo_r} respectively are bijective for all $n \in \IZ$.  The original
statement of the Baum-Connes Conjecture with trivial coefficients can be found
in~\cite[Conjecture~3.15 on page~254]{Baum-Connes-Higson(1994)}.

Our formulation of these conjectures follows the homotopy theoretic approach in
\cite{Davis-Lueck(1998)}.  The original assembly maps are defined differently. We do not
give the proof that our maps agree with the original ones but at least refer to
\cite[page~239]{Davis-Lueck(1998)}, where the Farrell-Jones Conjecture is treated and to
Hambleton-Pedersen~\cite{Hambleton-Pedersen(2004)}, where such identification is given for
the Baum-Connes Conjecture with trivial coefficients. 


\subsection{Inheritance under colimits}
\label{subsec:Inheritance_under_colimits}

The main purpose of this paper is to prove that these conjectures are inherited under
colimits over directed systems of groups (with not necessarily injective structure maps).
We want to show:

\begin{theorem}[Inheritance under colimits]
\label{the:inheritance_under_colimits}
Let $\{G_i \mid i \in I\}$ be a directed system of groups with (not necessarily injective)
structure maps $\phi_{i,j} \colon G_i \to G_j$.  Let $G = \colim_{i \in I} G_i$ be its
colimit with structure maps $\psi_i \colon G_i \to G$. Let $R$ be a ring (with involution)
and let $A$ be a $C^*$-algebra with structure preserving $G$-action. Given $i \in I$ and a
subgroup $H \subseteq G_i$, we let $H$ act on $R$ and $A$ by restriction with the group
homomorphism $(\psi_i)|_H \colon H \to G$. Fix $n \in \IZ$. Then:

\begin{enumerate}

\item \label{the:inheritance_under_colimits:arbitrary} If the assembly map
  $$\asmb^H_n \colon H_n^H(\EGF{H}{\VCyc};\bfK_R) \to H_n^H(\pt;\bfK_R) = K_n(R\rtimes
  H)$$
  of~\eqref{ass_K_n} is bijective for all $n \in \IZ$, all $i \in I$ and all
  subgroups $H \subseteq G_i$, then for every subgroup $K \subseteq G$ of $G$ the assembly
  map
  $$\asmb^K_n \colon H_n^K(\EGF{K}{\VCyc};\bfK_R) \to H_n^K(\pt;\bfK_R) = K_n(R\rtimes
  K)$$
  of~\eqref{ass_K_n} is bijective for all $n \in \IZ$.

  The corresponding version is true for the assembly maps given in~\eqref{ass_KH_n},
  \eqref{ass_L_n}, \eqref{ass_K_ntopo_l1}, and~\eqref{ass_K_ntopo_m};

\item \label{the:inheritance_under_colimits:injective} Suppose that all structure maps
  $\phi_{i,j}$ are injective and that the assembly map
  $$\asmb^{G_i}_n \colon H_n^{G_i}(\EGF{G_i}{\VCyc};\bfK_R) \to H_n^{G_i}(\pt;\bfK_R) =
  K_n(R\rtimes G_i)$$
  of~\eqref{ass_K_n} is bijective for all $n \in \IZ$ and $i \in I$.
  Then the assembly map
  $$\asmb^G_n \colon H_n^G(\EGF{G}{\VCyc};\bfK_R) \to H_n^G(\pt;\bfK_R) = K_n(R\rtimes
  G)$$
  of~\eqref{ass_K_n} is bijective for all $n \in \IZ$;

  The corresponding statement is true for the assembly maps given in~\eqref{ass_KH_n},
  \eqref{ass_L_n}, \eqref{ass_K_ntopo_l1}, \eqref{ass_K_ntopo_r},
  and~\eqref{ass_K_ntopo_m}.
\end{enumerate}
\end{theorem}

Theorem~\ref{the:inheritance_under_colimits} will follow from
Theorem~\ref{the:isomorphism_conjecture_is_stable_under_colim} and
Lemma~\ref{lem:homology_theories_are_(strongly)_continuous} as soon as we have proved
Theorem~\ref{the:Construction_of_equivariant_homology_theories}. Notice that the
version~\eqref{ass_K_ntopo_r} does not appear in
assertion~\ref{the:inheritance_under_colimits:arbitrary}. A counterexample will be
discussed below.  The (fibered) version of
Theorem~\ref{the:inheritance_under_colimits}~\ref{the:inheritance_under_colimits:arbitrary}
in the case of algebraic $K$-theory and $L$-theory with coefficients in $\IZ$ with trivial
$G$-action has been proved by Farrell-Linnell~\cite[Theorem~7.1]{Farrell-Linnell(2003a)}.


\subsection{Colimits of hyperbolic groups}
\label{subsec:Colimits_of_hyperbolic_groups}

In~\cite[Section~7]{Higson-Lafforgue-Skandalis(2002)} Higson, Lafforgue and Skandalis
construct counterexamples to the \emph{Baum-Connes-Conjecture with coefficients}, actually
with a commutative $C^*$-algebra as coefficients. 
They formulate precise properties for 
a group $G$ which imply that it does \emph{not} satisfy the Baum-Connes Conjecture with
coefficients.
Gromov~\cite{Gromov(2000)} describes the construction of such a group $G$
as a colimit over a directed system of groups $\{G_i \mid i \in I\}$, where each $G_i$ is
hyperbolic.

This construction did raise the hope that these groups $G$ may also be counterexamples to
the Baum-Connes Conjecture with trivial coefficients. But --- to the authors' knowledge
--- this has not been proved and no counterexample to the Baum-Connes Conjecture with
trivial coefficients is known.

Of course one may wonder whether such counterexamples to the Baum-Connes Conjecture with
coefficients or with trivial coefficients respectively may also be counterexamples to the
Farrell-Jones Conjecture or the Bost Conjecture with coefficients or with trivial
coefficients respectively. However, this can be excluded by the following result.

\begin{theorem} \label{the:Conjectures_for_colimits_of_hyperbolic_groups}
  Let $G$ be the colimit of the directed system $\{G_i \mid i \in I\}$ of groups (with not
  necessarily injective structure maps). Suppose that each $G_i$ is hyperbolic.
  Let $K \subseteq G$ be a subgroup. Then:

\begin{enumerate}
\item \label{the:Conjectures_for_colimits_of_hyperbolic_groups:K} The group $K$ satisfies
  for every ring $R$ on which $K$ acts by ring automorphisms the Farrell-Jones Conjecture
  for algebraic $K$-theory with coefficients in $R$, i.e., the assembly
  map~\eqref{ass_K_n} is bijective for all $n \in \IZ$;

\item \label{the:Conjectures_for_colimits_of_hyperbolic_groups:KH} The group $K$ satisfies
  for every ring $R$ on which $K$ acts by ring automorphisms the Farrell-Jones Conjecture
  for homotopy $K$-theory with coefficients in $R$, i.e., the assembly
  map~\eqref{ass_KH_n} is bijective for all $n \in \IZ$;

\item \label{the:Conjectures_for_colimits_of_hyperbolic_groups:l1} The group $K$ satisfies
  for every $C^*$-algebra $A$ on which $K$ acts by $C^*$-auto\-morphisms the Bost Conjecture
  with coefficients in $A$, i.e., the assembly map~\eqref{ass_K_ntopo_l1} is bijective for
  all $n \in \IZ$.

 \end{enumerate}

\end{theorem}

\begin{proof}
If $G$ is the colimit of the directed system $\{G_i \mid i \in
I\}$, then the subgroup $K \subseteq G$ is the colimit of the
directed system $\{\psi_i^{-1}(K) \mid i \in I\}$, where $\psi_i
\colon G_i \to G$ is the structure map.  Hence it suffices to
prove
Theorem~\ref{the:Conjectures_for_colimits_of_hyperbolic_groups} in
the case $G = K$. This case follows from
Theorem~\ref{the:inheritance_under_colimits}~\ref{the:inheritance_under_colimits:arbitrary}
as soon as one can show that the Farrell-Jones Conjecture for
algebraic $K$-theory, the Farrell-Jones Conjecture for homotopy
$K$-theory, or the Bost Conjecture respectively holds for every
subgroup $H$ of a hyperbolic group $G$ with arbitrary coefficients
$R$ and $A$ respectively.

Firstly we prove this for the Bost Conjecture.  Mineyev and
Yu~\cite[Theorem~17]{Mineyev-Yu(2002)} show that every hyperbolic group G admits a
G-invariant metric $\widehat{d}$ which is weakly geodesic and strongly bolic. Since every
subgroup $H$ of $G$ clearly acts properly on G with respect to any discrete metric, it
follows that $H$ belongs to the class $\calc'$ as described by Lafforgue
in~\cite[page~5]{Lafforgue(2002)} (see also the remarks at the top of page~6
in~\cite{Lafforgue(2002)}).  
Now the claim is a direct consequence of~\cite[Theorem~0.0.2]{Lafforgue(2002)}.

The claim for the Farrell-Jones Conjecture is proved for algebraic $K$-theory and homotopy
$K$-theory in Bartels-L\"uck-Reich~\cite{Bartels-Lueck-Reich(2007appl)} which is based on
the results of~\cite{Bartels-Lueck-Reich(2007hyper)}.
\end{proof}

There are further groups with unusual properties that can be obtained as
colimits of hyperbolic groups.
This class contains for instance a torsion-free non-cyclic group all whose proper subgroups
are cyclic constructed by Ol'shanskii~\cite{Olshanskii(1979)}.
Further examples are mentioned in \cite[p.5]{Olshanskii-Osin-Sapir(2007)}
and \cite[Section~4]{Sapir(2007)}. 

We mention that if one can prove the $L$-theoretic version of the Farrell-Jones Conjecture
for subgroups of hyperbolic groups with arbitrary coefficients, then it is also true for
subgroups of colimits of hyperbolic groups by the argument above.


\subsection{Discussion of (potential) counterexamples}
\label{subsec:Discussion_of_(potential)_counterexamples}

If $G$ is an infinite group which satisfies Kazhdan's property (T), then the assembly
map~\eqref{ass_K_ntopo_m} for the maximal group $C^*$-algebra fails to be an isomorphism
if the assembly map~\eqref{ass_K_ntopo_r} for the reduced group $C^*$-algebra is injective
(which is true for a very large class of groups and in particular for all hyperbolic
groups by~\cite{Kasparov-Skandalis(2003)}).  The reason is that a group has property (T)
if and only if the trivial representation $1_G$ is isolated in the dual $\widehat{G}$ of
$G$.  This implies that $C^*_m(G)$ has a splitting $\IC\oplus \ker(1_G)$, where we regard
$1_G$ as a representation of $C^*_m(G)$.  If $G$ is infinite, then the first summand is in
the kernel of the regular representation $\lambda\colon C^*_m(G)\to C_r^*(G)$ (see for
instance~\cite{Harpe-Valette(1989)}), hence the $K$-theory map $\lambda \colon
K_0(C^*_m(G))\to K_0(C_r^*(G))$ is not injective. 
For a finite group $H$ we have $A \rtimes_r H = A \rtimes_m H$
and hence we can apply \cite[Lemma~4.6]{Davis-Lueck(1998)} to
identify the domains of~\eqref{ass_K_ntopo_r} and~\eqref{ass_K_ntopo_m}.
Under this identification  
the composition of the assembly
map~\eqref{ass_K_ntopo_m} with $\lambda$ is the assembly map~\eqref{ass_K_ntopo_r} and the
claim follows.

Hence the Baum-Connes Conjecture for the maximal group $C^*$-algebras is not true in
general since the Baum-Connes Conjecture for the reduced group
$C^*$-algebras has been proved for some
groups with property (T) by Lafforgue~\cite{Lafforgue(1998)} (see also \cite{Skandalis(1999)}).
So in the sequel our discussion refers always to the Baum-Connes Conjecture for
the reduced group $C^*$-algebra.

One may speculate that the Baum-Connes Conjecture with trivial coefficients is less likely
to be true for a given group $G$ than the Farrell-Jones Conjecture or the Bost Conjecture.
Some evidence for this speculation comes from lack of functoriality of the reduced group
$C^*$-algebra. A group homomorphism $\alpha \colon H \to G$ induces in general not a
$C^*$-homomorphism $C^*_r(H) \to C^*_r(G)$, one has to require that its kernel is amenable.
Here is a counterexample, namely, if $F$ is a non-abelian free group, then $C^*_r(F)$ is a
simple algebra \cite{Powers(1975)} and hence there is no unital algebra homomorphism
$C^*_r(F) \to C^{\ast}_r ( \{ 1 \})= \IC $.  On the other hand, any group homomorphism
$\alpha \colon H \to G$ induces a homomorphism
$$H_n^H(\EGF{H}{\Fin};\bfK^{\topo}_{\IC,r}) \xrightarrow{\ind_{\alpha}}
H_n^G(\alpha_*\EGF{H}{\Fin};\bfK^{\topo}_{\IC,r}) \xrightarrow{H_n^G(f)}
H_n^G(\EGF{G}{\Fin};\bfK^{\topo}_{\IC,r})$$
where $G$ acts trivially on $\IC$ and $f
\colon \alpha_* \EGF{H}{\Fin} \to \EGF{G}{\Fin}$ is the up to $G$-homotopy unique $G$-map.
Notice that the induction map $\ind_{\alpha}$ exists since the isotropy groups of
$\EGF{H}{\Fin}$ are finite.  Moreover, this map is compatible under the assembly maps for
$H$ and $G$ with the map $K_n(C^*_r(\alpha)) \colon K_n(C^*_r(H)) \to K_n(C^*_r(G))$
provided that $\alpha$ has amenable kernel 
and hence $C^*_r(\alpha)$ is defined.  So the
Baum-Connes Conjecture implies that every group homomorphism $\alpha \colon H \to G$
induces a group homomorphism $\alpha_* \colon K_n(C^*_r(H)) \to K_n(C^*_r(G))$, although
there may be no $C^*$-homomorphism $C^*_r(H) \to C^*_r(G)$ induced by $\alpha$. No such
direct construction of $\alpha_*$ is known in general.

Here is another failure of the reduced group $C^*$-algebra.  Let $G$
be the colimit of the directed system $\{G_i \mid i \in I\}$ of groups
(with not necessarily injective structure maps). Suppose that for
every $i \in I$ and preimage $H$ of a finite group under the canonical map 
$\psi_i \colon G_i \to G$ the
Baum-Connes Conjecture for the maximal group $C^*$-algebra holds (This
is for instance true by~\cite{Higson-Kasparov(2001)} if $\ker(\psi_i)$
has the Haagerup property). Then
\begin{multline*}
  \colim_{i \in I} H_n^{G_i}(\EGF{G_i}{\Fin};\bfK^{\topo}_{\IC,m}) \xrightarrow{\cong}
  \colim_{i \in I} H_n^{G_i}(\EGF{G_i}{\psi_i^*\Fin};\bfK^{\topo}_{\IC,m})
  \\
  \xrightarrow{\cong} H_n^G(\EGF{G}{\Fin};\bfK^{\topo}_{\IC,m})
\end{multline*}
is a composition of two isomorphisms. The first map is bijective by
the Transitivity Principle~\ref{the:transitivity}, 
the second by 
Lemma~\ref{lem:directed_limits_and_calh?_ast} and 
Lemma~\ref{lem:homology_theories_are_(strongly)_continuous}.  
This implies that the following composition is an isomorphism
\begin{multline*}
  \colim_{i \in I} H_n^{G_i}(\EGF{G_i}{\Fin};\bfK^{\topo}_{\IC,r}) \to
  \colim_{i \in I} H_n^{G_i}(\EGF{G_i}{\psi_i^*\Fin};\bfK^{\topo}_{\IC,r})
  \\
  \to H_n^G(\EGF{G}{\Fin};\bfK^{\topo}_{\IC,r})
\end{multline*}
Namely, these two compositions are compatible with the passage from
the maximal to the reduced setting. This passage induces on the source
and on the target isomorphisms since $\EGF{G_i}{\Fin}$ and
$\EGF{G}{\Fin}$ have finite isotropy groups, for a finite group $H$ we
have $C^*_r(H) = C^*_m(H)$ and hence we can apply
\cite[Lemma~4.6]{Davis-Lueck(1998)}.  Now assume furthermore that the
Baum-Connes Conjecture for the reduced group $C^*$-algebra holds for
$G_i$ for each $i \in I$ and for $G$.  Then we obtain an isomorphism
$$\colim_{i \in I} K_n(C^*_r(G_i)) \xrightarrow{\cong}
K_n(C^*_r(G)).$$
Again it is in general not at all clear whether there
exists such a map in the case, where the structure maps $\psi_i \colon
G_i \to G$ do not have amenable kernels and hence do not induce maps
$C^*_r(G_i) \to C^*_r(G)$.

These arguments do not apply to the Farrell-Jones Conjecture or the Bost Conjecture.
Namely any group homomorphism $\alpha \colon H \to G$ induces maps $R \rtimes H \to R
\rtimes G$, $A\rtimes_{l^1} H \to A \rtimes_{l^1} G$, and $A\rtimes_m H \to A \rtimes_m G$
for a ring $R$ or a $C^*$-algebra $A$ with structure preserving $G$-action, where we equip
$R$ and $A$ with the $H$-action coming from $\alpha$. Moreover we will show for a directed
system $\{G_i \mid i \in I\}$ of groups (with not necessarily injective structure maps)
and $G = \colim_{i \in I} G_i$ that there are canonical isomorphisms
(see~Lemma~\ref{lem:homology_theories_are_(strongly)_continuous})
$$\begin{array}{lcl} \colim_{i \in I} K_n(R\rtimes G_i) & \xrightarrow{\cong} &
  K_n(R\rtimes G);
  \\
  \colim_{i \in I} \KH_n(R\rtimes G_i) & \xrightarrow{\cong} & \KH_n(R\rtimes G);
  \\
  \colim_{i \in I} L_n^{\langle - \infty \rangle}(R\rtimes G_i) & \xrightarrow{\cong} &
  L_n^{\langle - \infty \rangle}(R\rtimes G);
  \\
  \colim_{i \in I} K_n(A \rtimes_{l^1} G_i) & \xrightarrow{\cong} & K_n(A \rtimes_{l^1}
  G);
  \\
  \colim_{i \in I} K_n(A \rtimes_{m} G_i) & \xrightarrow{\cong} & K_n(A \rtimes_{m} G).
\end{array}$$

Let $A$ be a $C^*$-algebra with $G$-action by $C^*$-automorphisms.  We can consider $A$ as
a ring only. Notice that we get a commutative diagram

$$\xymatrix{H_n^G(\EGF{G}{\VCyc};\bfK_A) \ar[d] \ar[r] & \KH_n(A\rtimes G) \ar[d]
  \\
  H_n^G(\EGF{G}{\VCyc};\bfKH_A) \ar[r] & \KH_n(A\rtimes G)
  \\
  H_n^G(\EGF{G}{\Fin};\bfKH_A) \ar[d] \ar[u]_{\cong} \ar[r] & \KH_n(A\rtimes G) \ar[d]
  \ar[u]_{\id}
  \\
  H_n^G(\EGF{G}{\Fin};\bfK^{\topo}_{A,l^1}) \ar[d]_{\cong} \ar[r] & K_n(A \rtimes_{l^1} G)
  \ar[d]
  \\
  H_n^G(\EGF{G}{\Fin};\bfK^{\topo}_{A,m}) \ar[d]_{\cong} \ar[r] & K_n(A \rtimes_m G)\ar[d]
  \\
  H_n^G(\EGF{G}{\Fin};\bfK^{\topo}_{A,r}) \ar[r] & K_n(A \rtimes_r G) }$$
where the
horizontal maps are assembly maps and the vertical maps are change of theory and rings
maps or induced by the up to $G$-homotopy unique $G$-map $\EGF{G}{\Fin} \to
\EGF{G}{\VCyc}$.  The second left vertical map, which is marked with $\cong$, is
bijective.  This is shown in \cite[Remark~7.4]{Bartels-Lueck(2006)} in the
case, where $G$ acts trivially on $R$, the proof carries directly over to the general
case.
The fourth
and fifth vertical left arrow, which are marked with $\cong$, are bijective, since for a
finite group $H$ we have $A \rtimes H = A\rtimes_{l^1} H = A \rtimes_r H = A \rtimes_m H$
and hence we can apply \cite[Lemma~4.6]{Davis-Lueck(1998)}. 
In particular the Bost
Conjecture and the Baum-Connes Conjecture together imply that the map $K_n(A \rtimes_{l^1}
G) \to K_n(A \rtimes_r G)$ is bijective, the map $K_n(A \rtimes_{l^1} G) \to K_n(A
\rtimes_m G)$ is split injective and the map $K_n(A \rtimes_m G) \to K_n(A \rtimes_r G)$
is split surjective.

The upshot of this discussions is:

\begin{itemize}

\item The counterexamples of Higson, Lafforgue and
  Skandalis~\cite[Section~7]{Higson-Lafforgue-Skandalis(2002)} to the Baum-Connes
  Conjecture with coefficients are not counterexamples to the Farrell-Jones Conjecture or
  the Bost Conjecture;

\item The counterexamples of Higson, Lafforgue and
  Skandalis~\cite[Section~7]{Higson-Lafforgue-Skandalis(2002)} show that the map $K_n(A
  \rtimes_{l^1} G) \to K_n(A \rtimes_r G)$ is in general not bijective;

\item The passage from the topological $K$-theory of the Banach algebra $l^1(G)$ to the
  reduced group $C^*$-algebra is problematic and may cause failures of the Baum-Connes
  Conjecture;

\item The Bost Conjecture and the Farrell-Jones Conjecture are more likely to be true than
  the Baum-Connes Conjecture;

\item There is --- to the authors' knowledge --- no promising candidate of a group $G$ for
  which a strategy is in sight to show that the Farrell-Jones Conjecture or the Bost
  Conjecture are false. (Whether it is reasonable to believe that these conjectures are
  true for all groups is a different question.)

\end{itemize}


\subsection{Homology theories and spectra}
\label{subsec:Homology_theories_and_spectra}

The general strategy of this paper is to present most of the
arguments in terms of equivariant homology theories. Many of the
arguments for the Farrell-Jones Conjecture, the Bost Conjecture or
the Baum-Connes Conjecture become the same, the only difference
lies in the homology theory we apply them to. This is convenient
for a reader who is not so familiar with spectra and prefers to
think of $K$-groups and not of $K$-spectra.

The construction of these equivariant homology theories is a second
step and done in terms of spectra.  Spectra cannot be avoided in
algebraic $K$-theory by definition and since we want to compare also
algebraic and topological $K$-theory, we need spectra descriptions
here as well. Another nice feature of the approach to equivariant
topological $K$-theory via spectra is that it yields a theory which
can be applied to all $G$-$CW$-complexes.
This will allow us to consider
in the case $G = \colim_{i \in I} G_i$ the equivariant
$K$-homology of the $G_i$-$CW$-complex 
$\psi^*_i\EGF{G}{\Fin} = \EGF{G_i}{\psi^*\Fin}$
although $\psi^*_i\EGF{G}{\Fin}$ has infinite isotropy groups
if the structure map $\psi_i \colon G_i \to G$ has infinite kernel.

Details of the constructions of the relevant spectra, namely, the proof of
Theorem~\ref{the:necessary_functors_from_groupoids_to_spectra}, will be deferred
to~\cite{Bartels-Joachim-Lueck(2007)}. We will use the existence of
these spectra as a black box. These constructions require some work
and technical skills, but their details are not at all relevant for
the results and ideas of this paper and their existence is not at all surprising.


\subsection{Twisting by cocycles}
\label{subsec:Twisting_by_cocycles}

In the $L$-theory case one encounters also non-orien\-table
manifolds. In this case twisting with the first Stiefel-Whitney class
is required. In a more general setup one is given a group $G$, a ring
$R$ with involution and a group homomorphism $w\colon G \to
\cent(R^{\times})$ to the center of the multiplicative group of units
in $R$. So far we have used the standard involution on the group ring
$RG$, which is given by $\overline{r\cdot g} = \overline{r} \cdot
g^{-1}$. One may also consider the $w$-twisted involution given by
$\overline{r\cdot g} = \overline{r}w(g) \cdot g^{-1}$.  All the
results in this paper generalize directly to this case since one can
construct a modified $L$-theory spectrum functor (over $G$) using the
$w$-twisted involution and then the homology arguments are just
applied to the equivariant homology theory associated to this
$w$-twisted $L$-theory spectrum.


\subsection{Acknowledgements}
\label{subsec:Acknowledgements}
The work was financially supported by the Sonderforschungsbereich 478 \--- Geometrische
Strukturen in der Mathematik \--- and the Max-Planck-Forschungspreis of the third author.



\typeout{------------------- Section 1: Equivariant homology theories ------------------}

\section{Equivariant homology theories}
\label{sec:Equivariant_homology_theories}

In this section we briefly explain basic axioms, notions and facts about equivariant
homology theories as needed for the purposes of this article. The main examples which will
play a role in connection with the Bost, the Baum-Connes and the Farrell-Jones Conjecture
will be presented later in Theorem~\ref{the:Construction_of_equivariant_homology_theories}.

Fix a group $G$ and a ring $\Lambda$. In most cases $\Lambda$ will be $\IZ$.  The following
definition is taken from~\cite[Section~1]{Lueck(2002b)}.

\begin{definition}[$G$-homology theory] \label{def:G-homology-theory}
  A \emph{$G$-homology theory $\calh_*^G$ with values in $\Lambda$-modules} is a
  collection of covariant functors $\calh^G_n$ from the category of $G$-$CW$-pairs to the
  category of $\Lambda$-modules indexed by $n \in \IZ$ together with natural
  transformations $\partial_n^G(X,A)\colon \calh_n^G(X,A) \to \calh_{n-1}^G(A):=
  \calh_{n-1}^G(A,\emptyset)$ for $n \in \IZ$ such that the following axioms are
  satisfied:
\begin{itemize}

\item $G$-homotopy invariance \\[1mm]
  If $f_0$ and $f_1$ are $G$-homotopic maps $(X,A) \to (Y,B)$ of $G$-$CW$-pairs, then
  $\calh_n^G(f_0) = \calh^G_n(f_1)$ for $n \in \IZ$;

\item Long exact sequence of a pair\\[1mm]
  Given a pair $(X,A)$ of $G$-$CW$-complexes, there is a long exact sequence
  \begin{multline*} \quad \quad \quad \ldots \xrightarrow{\calh^G_{n+1}(j)} \calh_{n+1}^G(X,A)
    \xrightarrow{\partial_{n+1}^G} \calh_n^G(A) \xrightarrow{\calh^G_{n}(i)} \calh_n^G(X)
    \\
    \xrightarrow{\calh^G_{n}(j)} \calh_n^G(X,A) \xrightarrow{\partial_n^G} \ldots,
  \end{multline*}
  where $i\colon A \to X$ and $j\colon X \to (X,A)$ are the inclusions;

\item Excision \\[1mm]
  Let $(X,A)$ be a $G$-$CW$-pair and let $f\colon A \to B$ be a cellular $G$-map of
  $G$-$CW$-complexes. Equip $(X\cup_f B,B)$ with the induced structure of a $G$-$CW$-pair.
  Then the canonical map $(F,f)\colon (X,A) \to (X\cup_f B,B)$ induces an isomorphism
  $$\calh_n^G(F,f)\colon \calh_n^G(X,A) \xrightarrow{\cong} \calh_n^G(X\cup_f B,B);$$

\item Disjoint union axiom\\[1mm]
  Let $\{X_i\mid i \in I\}$ be a family of $G$-$CW$-complexes. Denote by $j_i\colon X_i
  \to \coprod_{i \in I} X_i$ the canonical inclusion.  Then the map
  $$\bigoplus_{i \in I} \calh^G_{n}(j_i)\colon \bigoplus_{i \in I} \calh_n^G(X_i)
  \xrightarrow{\cong} \calh_n^G\left(\coprod_{i \in I} X_i\right)$$
  is bijective.

\end{itemize}

\end{definition}

Let $\calh^G_*$ and $\calk^G_*$ be $G$-homology theories.  A \emph{natural
  transformation $T_* \colon \calh^G_* \to \calk^G_*$ of $G$-homology theories} is a
sequence of natural transformations $T_n \colon \calh_n^G \to \calk_n^G$ of functors from
the category of $G$-$CW$-pairs to the category of $\Lambda$-modules which are compatible
with the boundary homomorphisms.

\begin{lemma} \label{lem:transf_between_G-homologies}
  Let $T_* \colon \calh^G_* \to \calk^G_*$ be a natural transformation of $G$-homology theories.
  Suppose that $T_n(G/H)$ is bijective for every homogeneous space $G/H$ and $n \in \IZ$.

  Then $T_n(X,A)$ is bijective for every $G$-$CW$-pair $(X,A)$ and $n \in \IZ$.
\end{lemma}
\begin{proof}
  The disjoint union axiom implies that both $G$-homology theories are compatible with
  colimits over directed systems indexed by the natural numbers (such as the system given
  by the skeletal filtration $X_0 \subseteq X_1 \subseteq X_2 \subseteq \ldots \subseteq
  \cup_{n \ge 0} X_n = X$).  The argument for this claim is analogous to the one in
  \cite[7.53]{Switzer(1975)}.  Hence it suffices to prove the bijectivity for
  finite-dimensional pairs. Using the axioms of a $G$-homology theory, the five lemma and
  induction over the dimension one reduces the proof to the special case $(X,A) =
  (G/H,\emptyset)$.
\end{proof}

Next we present a slight variation of the notion of an equivariant homology theory
introduced in~\cite[Section~1]{Lueck(2002b)}. We have to treat this variation since we
later want to study coefficients over a fixed group $\Gamma$ which we will then pullback via
group homomorphisms with $\Gamma$ as target.  Namely, fix a group $\Gamma$.  A group
$(G,\xi)$ over $\Gamma$ is a group $G$ together with a group homomorphism $\xi \colon G
\to \Gamma$.  A map $\alpha \colon (G_1,\xi_1) \to (G_2,\xi_2)$ of groups over $\Gamma$ is
a group homomorphisms $\alpha \colon G_1 \to G_2$ satisfying $\xi_2 \circ \alpha = \xi_1$.

Let $\alpha\colon H \to G$ be a group homomorphism.  Given an $H$-space $X$, define the
\emph{induction of $X$ with $\alpha$} to be the $G$-space denoted by $\alpha_* X$ which is
the quotient of $G \times X$ by the right $H$-action $(g,x) \cdot h := (g\alpha(h),h^{-1}
x)$ for $h \in H$ and $(g,x) \in G \times X$.  If $\alpha\colon H \to G$ is an inclusion,
we also write $\ind_H^G$ instead of $\alpha_*$. If $(X,A)$ is an $H$-$CW$-pair, then
$\alpha_*(X,A)$ is a $G$-$CW$-pair.

\begin{definition}[Equivariant homology theory over a group $\Gamma$]
\label{def:Equivariant_homology_theory_over_a_group_Gamma}
An \emph{equivariant homology theory $\calh^?_*$ with values in $\Lambda$-modules over a
  group $\Gamma$} assigns to every group $(G,\xi)$ over $\Gamma$ a $G$-homology theory
$\calh^G_*$ with values in $\Lambda$-modules and comes with the following so called
\emph{induction structure}: given a homomorphism $\alpha \colon (H,\xi) \to (G,\mu)$ of
groups over $\Gamma$ and an $H$-$CW$-pair $(X,A)$, there are for each $n \in \IZ$ natural
homomorphisms
\begin{eqnarray}
&\ind_{\alpha}\colon \calh_n^H(X,A)
\to
\calh_n^G(\alpha_*(X,A))&
\label{induction_structure}
\end{eqnarray}
satisfying

\begin{itemize}

\item Compatibility with the boundary homomorphisms\\[1mm]
  $\partial_n^G \circ \ind_{\alpha} = \ind_{\alpha} \circ \partial_n^H$;

\item Functoriality\\[1mm]
  Let $\beta \colon (G,\mu) \to (K,\nu)$ be another morphism of groups over $\Gamma$.
  Then we have for $n \in \IZ$
  $$\ind_{\beta \circ \alpha} = \calh^K_n(f_1)\circ\ind_{\beta} \circ \ind_{\alpha}\colon
  \calh^H_n(X,A) \to \calh_n^K((\beta\circ\alpha)_*(X,A)),$$
  where $f_1\colon
  \beta_*\alpha_*(X,A) \xrightarrow{\cong} (\beta\circ \alpha)_*(X,A), \quad (k,g,x)
  \mapsto (k\beta(g),x)$ is the natural $K$-homeo\-mor\-phism;

\item Compatibility with conjugation\\[1mm]
  Let $(G,\xi)$ be a group over $\Gamma$ and let $g \in G$ be an element with $\xi(g) =
  1$.  Then the conjugation homomorphisms $c(g) \colon G \to G$ defines a morphism $c(g)
  \colon (G,\xi) \to (G,\xi)$ of groups over $\Gamma$.  Let $f_2\colon (X,A) \to
  c(g)_*(X,A)$ be the $G$-homeomorphism which sends $x$ to $(1,g^{-1}x)$ in
  $G\times_{c(g)} (X,A)$.

  Then for every $n \in \IZ$ and every $G$-$CW$-pair $(X,A)$ the homomorphism
  $\ind_{c(g)}\colon \calh^G_n(X,A)\to \calh^G_n(c(g)_*(X,A))$ agrees with
  $\calh_n^G(f_2)$.

\item Bijectivity\\[1mm]
  If $\alpha \colon (H,\xi) \to (G,\mu)$ is a morphism of groups over $\Gamma$ such that
  the underlying group homomorphism $\alpha \colon H \to G$ is an inclusion of groups,
  then $\ind_{\alpha}\colon \calh_n^H(\pt) \to \calh_n^G(\alpha_*\pt) = \calh_n^G(G/H)$ is
  bijective for all $n \in \IZ$.
\end{itemize}

\end{definition}

Definition~\ref{def:Equivariant_homology_theory_over_a_group_Gamma} reduces to the one of
an equivariant homology in~\cite[Section~1]{Lueck(2002b)} if one puts $\Gamma = \{1\}$.

\begin{lemma} \label{lem:induction_bijective}
  Let $\alpha \colon (H,\xi) \to (G,\mu)$ be a morphism of groups over $\Gamma$.  Let
  $(X,A)$ be an  $H$-$CW$-pair such that $\ker(\alpha)$ acts freely on $X-A$.  Then
\begin{eqnarray*}
&\ind_{\alpha}\colon \calh_n^H(X,A)
\to
\calh_n^G(\alpha_*(X,A))&
\end{eqnarray*}
is bijective for all $n \in \IZ$.
\end{lemma}
\begin{proof}
  Let $\calf$ be the set of all subgroups of $H$ whose intersection with $\ker(\alpha)$ is
  trivial.  Obviously, this is a family, i.e., closed under conjugation and taking
  subgroups.  A $H$-$CW$-pair $(X,A)$ is called a $\calf$-$H$-$CW$-pair if the isotropy group of
  any point in $X-A$ belongs to $\calf$. A $H$-$CW$-pair $(X,A)$ is a
  $\calf$-$H$-$CW$-pair if and only if $\ker(\alpha)$ acts freely on $X-A$.

  The $n$-skeleton of $\alpha_*(X,A)$ is $\alpha_*$ applied to the $n$-skeleton of
  $(X,A)$.  Let $(X,A)$ be an $H$-$CW$-pair and let $f\colon A \to B$ be a cellular $H$-map
  of $H$-$CW$-complexes.  Equip $(X\cup_f B,B)$ with the induced structure of a
  $H$-$CW$-pair.  Then there is an obvious natural isomorphism of $G$-$CW$-pairs
  $$\alpha_*(X\cup_f B,B) \xrightarrow{\cong} (\alpha_*X\cup_{\alpha_*f}
  \alpha_*B,\alpha_*B).$$

  Now we proceed as in the proof of Lemma~\ref{lem:transf_between_G-homologies} but now
  considering the transformations
\begin{eqnarray*}
&\ind_{\alpha}\colon \calh_n^H(X,A) \to \calh_n^G(\alpha_*(X,A))&
\end{eqnarray*}
only for $\calf$-$H$-$CW$-pairs $(X,A)$. Thus we can reduce the claim to the special case
$(X,A) = H/L$ for some subgroup $L \subseteq H$ with $L \cap \ker(\alpha) = \{1\}$.  This
special case follows from the following commutative diagram whose vertical arrows are
bijective by the axioms and whose upper horizontal arrow is bijective since $\alpha$
induces an isomorphism $\alpha|_L \colon L \to \alpha(L)$.
$$\xycomsquare{\calh_n^L(\pt)}{\ind_{\alpha|_L \colon L \to
    \alpha(L)}}{\calh_n^{\alpha(L)}(\pt)} {\ind_L^H}{\ind_{\alpha(L)}^G}
{\calh^H_n(H/L)}{\ind_{\alpha}}{\calh_n^G(\alpha_*H/L) = \calh^G_n(G/\alpha(L))}$$

\end{proof}


\typeout{------------ Section 2: Equivariant homology theories and colimits -----------}

\section{Equivariant homology theories and colimits}
\label{sec:Equivariant_homology_theories_and_colimits}

Fix a group $\Gamma$ and an equivariant homology theory $\calh^?_*$ with values in
$\Lambda$-modules over $\Gamma$.

Let $X$ be a $G$-$CW$-complex. Let $\alpha \colon H \to G$ be a group homomorphism.
Denote by $\alpha^*X$ the $H$-$CW$-complex obtained from $X$ by \emph{restriction with
  $\alpha$}. We have already introduced the induction $\alpha_*Y$ of an $H$-$CW$-complex
$Y$.  The functors ${\alpha}_*$ and $\alpha^*$ are adjoint to one another.  In particular
the adjoint of the identity on $\alpha^* X$ is a natural $G$-map
\begin{eqnarray}
& f(X,\alpha) \colon {\alpha}_*\alpha^*X \to X. &
\label{equ:f(X,alpha)}
\end{eqnarray}
It sends an element in $G \times_{\alpha} \alpha^*X$ given by $(g,x)$ to $gx$.

Consider a map $\alpha \colon (H,\xi) \to (G,\mu)$ of groups over $\Gamma$.  Define the
$\Lambda$-map
\begin{eqnarray*}
& a_n = a_n(X,\alpha)\colon \calh_n^{H}(\alpha^*X)
\xrightarrow{\ind_{\alpha}} \calh_n^G({\alpha}_*\alpha^*X)
\xrightarrow{\calh_n^G(f(X,\alpha))}
\calh_n^G(X).&
\end{eqnarray*}
If $\beta \colon (G,\mu) \to (K,\nu)$ is another morphism of groups over $\Gamma$, then by
the axioms of an induction structure the composite $\calh_n^{H}(\alpha^*\beta^*X)
\xrightarrow{a_n(\beta^*X,\alpha)} \calh_n^{G}(\beta^*X) \xrightarrow{a_n(X,\beta)}
\calh_n^{K}(X)$ agrees with $a_n(X,\beta \circ \alpha) \colon
\calh_n^{H}(\alpha^*\beta^*X) = \calh_n^{H}((\beta\circ \alpha)^*X) \to \calh_n^{G}(X)$
for a $K$-$CW$-complex $X$.

Consider a directed system of groups $\{G_i \mid i \in I\}$ with $G = \colim_{i \in I}
G_i$ and structure maps $\psi_i \colon G_i \to G$ for $i \in I$ and $\phi_{i,j} \colon G_i
\to G_j$ for $i,j \in I, i \le j$.  We obtain for every $G$-$CW$-complex $X$ a system of
$\Lambda$-modules $\{\calh^{G_i}(\psi_i^*X) \mid i \in I\}$ with structure maps
$a_n(\psi_j^*X,\phi_{i,j}) \colon \calh^{G_i}(\psi_i^*X) \to \calh^{G_j}(\psi_j^*X)$.  We
get a map of $\Lambda$-modules
\begin{eqnarray}
& & t_n^G(X,A)~:=~\colim_{i \in I} a_n(X,\psi_i) \colon \colim_{i \in
    I} \calh_n^{G_i}(\psi_i^*(X,A))~\to~\calh_n^G(X,A).
\label{t_nG(X)}
\end{eqnarray}

The next definition is an extension of \cite[Definition~3.1]{Bartels-Lueck(2006)}.

\begin{definition}[(Strongly) continuous equivariant homology theory]
\label{def:strongly_continuous_equivariant_homology_theory}
An equivariant homology theory $\calh^?_*$ over $\Gamma$ is called \emph{continuous} if
for every group $(G,\xi)$ over $\Gamma$ and every directed system of subgroups $\{G_i \mid
i \in I\}$ of $G$ with $G = \bigcup_{i \in I} G_i$ the $\Lambda$-map (see~\eqref{t_nG(X)})
$$t^G_n(\pt) \colon \colim_{i \in I} \calh^{G_i}_n(\pt) \to \calh^G_n(\pt)$$
is an
isomorphism for every $n \in \IZ$.

An equivariant homology theory $\calh^?_*$ over $\Gamma$ is called \emph{strongly
  continuous} if for every group $(G,\xi)$ over $\Gamma$ and every directed system of
groups $\{G_i \mid i \in I\}$ with $G = \colim_{i \in I} G_i$ and structure maps $\psi_i
\colon G_i \to G$ for $i \in I$ the $\Lambda$-map
$$t^G_n(\pt) \colon \colim_{i \in I} \calh^{G_i}_n(\pt) \to
\calh^G_n(\pt)$$
is an isomorphism for every $n \in \IZ$.
\end{definition}

Here and in the sequel we view $G_i$ as a group over $\Gamma$ by $\xi
\circ \psi_i \colon G_i \to \Gamma$ and $\psi_i \colon G_i \to G$ as a
morphism of groups over $\Gamma$.

\begin{lemma} \label{lem:directed_limits_and_calh?_ast}
  Let $(G,\xi)$ be a group over $\Gamma$. Consider a directed system
  of groups $\{G_i \mid i \in I\}$ with $G = \colim_{i \in I} G_i$ and
  structure maps $\psi_i \colon G_i \to G$ for $i \in I$.  Let $(X,A)$
  be a $G$-$CW$-pair.  Suppose that $\calh^?_*$ is strongly continuous.

  Then the $\Lambda$-homomorphism
  (see~\eqref{t_nG(X)})
  $$t_n^G(X,A) \colon \colim_{i \in I}
  \calh_n^{G_i}(\psi_i^*(X,A))~\xrightarrow{\cong}~\calh_n^G(X,A)$$
  is bijective for every $n \in \IZ$.
\end{lemma}
\begin{proof} The functor sending a directed systems of $\Lambda$-modules
  to its colimit is an exact functor and compatible with direct sums
  over arbitrary index maps. If $(X,A)$ is a pair of
  $G$-$CW$-complexes, then $(\psi_i^*X ,\psi_i^* A)$ is a pair of
  $G_i$-$CW$-complexes. Hence the collection of maps $\{t^G_n(X,A) \mid
  n \in \IZ\}$ is a natural transformation of $G$-homology theories of pairs
  of $G$-$CW$-complexes which satisfy the disjoint union axiom. Hence
  in order to show that $t^G_n(X,A)$ is bijective for all $n \in \IZ$
  and all pairs of $G$-$CW$-complexes $(X,A)$, it suffices by
  Lemma~\ref{lem:transf_between_G-homologies} to prove this
  in the special case $(X,A) = (G/H,\emptyset)$.

  For $i \in I$ let $k_i \colon G_i/\psi_i^{-1}(H) \to \psi_i^* (G/H)$
  be the $G_i$-map sending $g_i\psi_i^{-1}(H)$ to $\psi_i(g_i)H$.
  Consider a directed system of $\Lambda$-modules
  $\{\calh_n^{G_i}(G_i/\psi_i^{-1}(H))\mid i \in I\}$ whose structure
  maps for $i,j \in I, i\le j$ are given by the composite
  $$\calh_n^{G_i}(G_i/\psi_i^{-1}(H)) \xrightarrow{\ind_{\phi_{i,j}}}
  \calh_n^{G_j}(G_j \times_{\phi_{i,j}} G_i/\psi_i^{-1}(H))
  \xrightarrow{\calh_n^{G_j}(f_{i,j})}
  \calh_n^{G_j}(G_j/\psi_j^{-1}(H))$$
  for the $G_j$-map $f_{i,j}
  \colon G_j \times_{\phi_{i,j}} G_i/\psi_i^{-1}(H) \to
  G_j/\psi_j^{-1}(H)$ sending $(g_j,g_i\psi_i^{-1}(H))$ to
  $g_j\phi_{i,j}(g_i)\psi_j^{-1}(H)$. Then the following diagram
  commutes
  $$\xymatrix@C=25mm { \colim_{i \in I} \calh_n^{\psi_i^{-1}(H)}(\pt)
    \ar[r]^-{\colim_{i \in I} \ind_{\psi_i^{-1}(H)}^{G_i}}_{\cong}
    \ar[dd]_{t^H_n(\pt)}^{\cong} &
    \colim_{i \in I} \calh_n^{G_i}(G_i/\psi_i^{-1}(H)) \ar[d]^{\colim_{i \in I} \calh^{G_i}_n(k_i)} \\
    & \colim_{i \in I} \calh_n^{G_i}(\psi_i^*(G/H)) \ar[d]^{t^G_n(G/H)} \\
    \calh_n^H(\pt) \ar[r]^-{\ind_H^G}_{\cong} & \calh_n^G(G/H) }
  $$
  where the horizontal maps are the isomorphism given by induction.
  For the directed system $\{\psi_i^{-1}(H) \mid i \in I\}$ with
  structure maps $\phi_{i,j}|_{\psi^{-1}_i(H)} \colon \psi_i^{-1}(H)
  \to \psi^{-1}_j(H)$, the group homomorphism $\colim_{i \in I}
  \psi_i|_{\psi^{-1}_i(H)} \colon \colim_{i \in I} \psi^{-1}_i(H) \to
  H$ is an isomorphism. This follows by inspecting the standard model
  for the colimit over a directed system of groups. Hence the left
  vertical arrow is bijective since $\calh^?_*$ is strongly continuous
  by assumption. Therefore it remains to show that the map
\begin{equation}
\colim_{i \in I} \calh^{G_i}_n(k_i) \colon \colim_{i \in I} \calh_n^{G_i}(G_i/\psi_i^{-1}(H))
\to
 \colim_{i \in I} \calh_n^{G_i}(\psi_i^*G/H)
\label{colim_i_in_I_calhG_n(k_i)}
\end{equation}
is surjective.

Notice that the map given by the direct sum of the structure maps
$$
\bigoplus_{i \in I} \calh_n^{G_i}(\psi_i^*G/H) \to \colim_{i \in I}
\calh_n^{G_i}(\psi_i^*G/H)$$
is surjective. Hence it remains to show
for a fixed $i \in I$ that the image of the structure map
$$\calh_n^{G_i}(\psi_i^*G/H) \to \colim_{i \in I}
\calh_n^{G_i}(\psi_i^*G/H)$$
is contained in the image of the
map~\eqref{colim_i_in_I_calhG_n(k_i)}.

We have the decomposition of the $G_i$-set $\psi_i^*G/H$ into its
$G_i$-orbits
$$\coprod_{G_i(gH) \in G_i\backslash
  (\psi_i^*G/H)}~G_i/\psi_i^{-1}(gHg^{-1})
\xrightarrow{\cong}\psi_i^*G/H, \quad g_i\psi_i^{-1}(gHg^{-1}) \mapsto
\psi_i(g_i)gH.$$
It induces an identification of $\Lambda$-modules
$$
\bigoplus_{G_i(gH) \in G_i\backslash
  (\psi_i^*G/H)}~\calh^{G_i}_n\left(G_i/\psi_i^{-1}(gHg^{-1})\right)~=~\calh^{G_i}_n(\psi_i^*
G/H).$$
Hence it remains to show for fixed elements $i \in I$ and
$G_i(gH) \in G_i \backslash (\psi_i^*G/H)$ that the obvious
composition
$$\calh^{G_i}_n\left(G_i/\psi_i^{-1}(gHg^{-1})\right) \subseteq
\calh^{G_i}_n(\psi_i^* G/H) \to \colim_{i \in I}
\calh_n^{G_i}(\psi_i^*G/H)$$
is contained in the image of the
map~\eqref{colim_i_in_I_calhG_n(k_i)}.

Choose an index $j$ with $j \ge i$ and $g \in \im(\psi_j)$. Then the
structure map for $i \le j$ is a map $\calh^{G_i}_n(\psi_i^* G/H) \to
\calh^{G_j}_n(\psi_j^*G/H)$ which sends the summand corresponding to
$G_i(gH) \in G_i\backslash (\psi_i^*G/H)$ to the summand corresponding
to $G_j(1H) \in G_j\backslash (\psi_j^*G/H)$ which is by definition
the image of
$$\calh^{G_j}_n(k_j) \colon \calh^{G_j}_n(G_j/\psi_j^{-1}(H)) \to
\calh^{G_j}_n(\psi_j^*G/H).$$
Obviously the image of composite of the
last map with the structure map
$$\calh^{G_j}_n(\psi_j^*G/H) \to \colim_{i \in I}
\calh_n^{G_i}(\psi_i^*G/H)$$
is contained in the image of the
map~\eqref{colim_i_in_I_calhG_n(k_i)}. Hence the map
\eqref{colim_i_in_I_calhG_n(k_i)} is surjective. This finishes the
proof of Lemma~\ref{lem:directed_limits_and_calh?_ast}.
\end{proof}


\typeout{------------ Section 3: Isomorphism Conjectures and colimits
  ----------------}

\section{Isomorphism Conjectures and colimits}
\label{sec:Isomorphism_Conjectures_and_colimits}
A \emph{family $\calf$ of subgroups of $G$} is a collection of
subgroups of $G$ which is closed under conjugation and taking
subgroups.  Let $\EGF{G}{\calf}$ be the \emph{classifying space
  associated to $\calf$}.  It is uniquely characterized up to
$G$-homotopy by the properties that it is a $G$-$CW$-complex and that
$\EGF{G}{\calf}^H$ is contractible if $H \in \calf$ and is empty if $H
\notin \calf$.  For more information about these spaces
$\EGF{G}{\calf}$ we refer for instance to the survey
article~\cite{Lueck(2005s)}.  Given a group homomorphism $\phi \colon
K \to G$ and a family $\calf$ of subgroups of $G$, define the family
$\phi^*\calf$ of subgroups of $K$ by
\begin{eqnarray}
\phi^*\calf  & = & \{H \subseteq K \mid \phi(H) \in \calf\}.
\label{phiastcalf}
\end{eqnarray}
If $\phi$ is an inclusion of subgroups, we also write $\calf|_K$
instead of $\phi^*\calf$.

\begin{definition}[Isomorphism Conjecture for $\calh^?_*$]
\label{def:Isomorphism_Conjectures_for_calh?_ast}
Fix a group $\Gamma$ and an equivariant homology theory $\calh^?_*$
with values in $\Lambda$-modules over $\Gamma$.

A group $(G,\xi)$ over $\Gamma$ together with a family of subgroups
$\calf$ of $G$ satisfies the \emph{Isomorphism Conjecture (for
  $\calh^?_*$)} if the projection $\pr \colon \EGF{G}{\calf} \to \pt$
to the one-point-space $\pt$ induces an isomorphism
$$\calh^G_n(\pr) \colon \calh^G_n(\EGF{G}{\calf}) \xrightarrow{\cong}
\calh^G_n(\pt)$$
for all $n \in \IZ$.
\end{definition}

{}From now on fix a group $\Gamma$ and an equivariant homology theory
$\calh^?_*$ over $\Gamma$.

\begin{theorem}[Transitivity Principle] \label{the:transitivity}
  Let $(G,\xi)$ be a group over $\Gamma$.  Let $\calf \subseteq \calg$ be families of
  subgroups of $G$.  Assume that for every element $H \in \calg$ the group $(H,\xi|_H)$
  over $\Gamma$ satisfies the Isomorphism Conjecture for $\calf|_H$.

  Then the up to $G$-homotopy unique map 
  $\EGF{G}{\calf} \to \EGF{G}{\calg}$ induces an isomorphism 
  $\calh^G_n(\EGF{G}{\calf}) \to \calh^G_n(\EGF{G}{\calg})$
  for all $n \in \IZ$.
  In particular,
  $(G,\xi)$ satisfies the Isomorphism Conjecture for $\calg$ if and only if
  $(G,\xi)$ satisfies the Isomorphism Conjecture for $\calf$.
\end{theorem}

\begin{proof}
  The proof is completely analogous to the one
  in~\cite[Theorem~2.4, Lemma~2.2]{Bartels-Lueck(2006)}, where only the case
  $\Gamma = \{1\}$ is treated.
\end{proof}

\begin{theorem}
\label{the:IC_and_limits}
Let $(G,\xi)$ be a group over $\Gamma$. Let $\calf$ be a family of
subgroups of $G$.
\begin{enumerate}

\item \label{the:IC_and_limits:IC} Let $G$ be the directed union of
  subgroups $\{G_i \mid i \in I\}$. Suppose that $\calh^?_*$ is
  continuous and for every $i \in I$ the Isomorphism Conjecture holds
  for $(G_i,\xi|_{G_i})$ and $\calf|_{G_i}$.

  Then the Isomorphism Conjecture holds for $(G,\xi)$ and $\calf$;

\item \label{the:IC_and_limits:fibered_IC} Let $\{G_i \mid i \in I\}$ be a directed system
  of groups with $G = \colim_{i \in I} G_i$ and structure maps $\psi_i \colon G_i \to G$.
  Suppose that $\calh^?_*$ is strongly continuous and for every $i \in I$ the Isomorphism
  Conjecture holds for $(G_i,\xi \circ \psi_i)$ and $\psi_i^*\calf$.

  Then the Isomorphism Conjecture holds for $(G,\xi)$ and $\calf$.
\end{enumerate}
\end{theorem}
\begin{proof}\ref{the:IC_and_limits:IC}
  The proof is analogous to the one in
  \cite[Proposition~3.4]{Bartels-Lueck(2006)}.
  \\[1mm]\ref{the:IC_and_limits:fibered_IC} This follows from the
  following commutative square whose horizontal arrows are bijective
  because of Lemma~\ref{lem:directed_limits_and_calh?_ast} and the
  identification $\psi_i^*\EGF{G}{\calf} = \EGF{G_i}{\psi_i^*\calf}$
  $$\xymatrix@C=15mm {\colim_{i \in I}
    \calh_n^{G_i}(\EGF{G_i}{\psi_i^*\calf})
    \ar[r]^-{t_n^G(\EGF{G}{\calf})}_{\cong} \ar[d] &
    \calh_n^{G}(\EGF{G}{\calf}) \ar[d]
    \\
    \colim_{i \in I} \calh_n^{G_i}(\pt) \ar[r]^-{t_n^G(\pt)}_{\cong} &
    \calh_n^G(\pt) }
  $$
\end{proof}

Fix a class of groups $\calc$ closed under isomorphisms, taking
subgroups and taking quotients, e.g., the class of finite groups or
the class of virtually cyclic groups.  For a group $G$ let $\calc(G)$
be the family of subgroups of $G$ which belong to $\calc$.
\begin{theorem}
\label{the:isomorphism_conjecture_is_stable_under_colim}
Let $(G,\xi)$ be a group over $\Gamma$.
\begin{enumerate}
\item
\label{the:isomorphism_conjecture_is_stable_under_colim:injective}
Let $G$ be the directed union $G = \bigcup_{i \in I} G_i$ of subgroups
$G_i$ Suppose that $\calh^?_*$ is continuous and that the Isomorphism
Conjecture is true for $(G_i,\xi|_{G_i})$ and $\calc(G_i)$ for all $i \in I$.

Then the Isomorphism Conjecture is true for $(G,\xi)$ and $\calc(G)$;

\item
\label{the:isomorphism_conjecture_is_stable_under_colim:general}
Let $\{G_i \mid i \in I\}$ be a directed system of groups with $G =
\colim_{i \in I} G_i$ and structure maps $\psi_i \colon G_i \to G$.  
Suppose that $\calh^?_*$ is strongly
continuous and that the Isomorphism Conjecture is true for
$(H,\calc(H))$ for every $i \in I$ and every subgroup $H \subseteq
G_i$.

Then for every subgroup $K \subseteq G$ the Isomorphism Conjecture is true for 
$(K,\xi|_K)$ and $\calc(K)$.
\end{enumerate}
\end{theorem}
\begin{proof}\ref{the:isomorphism_conjecture_is_stable_under_colim:injective} 
  This follows from  Theorem~\ref{the:IC_and_limits}~\ref{the:IC_and_limits:IC} 
  since $\calc(G_i) = \calc(G)|_{G_i}$ holds for $i \in I$.
  \\[1mm]\ref{the:isomorphism_conjecture_is_stable_under_colim:general}
  If $G$ is the colimit of the directed system $\{G_i \mid i \in
  I\}$, then the subgroup $K \subseteq G$ is the colimit of the
  directed system $\{\psi_i^{-1}(K) \mid i \in I\}$.
  Hence we can assume $G = K$ without loss of generality.

  Since $\calc$ is closed under quotients by assumption, we have
  $\calc(G_i) \subseteq \psi_i^*\calc(G)$ for every $i \in I$.  Hence
  we can consider for any $i \in I$ the composition
  $$H_n^{G_i}(\EGF{G_i}{\calc(G_i)}) \to
  H_n^{G_i}(\EGF{G_i}{\psi_i^*\calc(G)}) \to H_n^{G_i}(\pt).$$
  Because
  of
  Theorem~\ref{the:IC_and_limits}~\ref{the:IC_and_limits:fibered_IC}
  it suffices to show that the second map is bijective. By assumption
  the composition of the two maps is bijective. Hence it remains to
  show that the first map is bijective. By
  Theorem~\ref{the:transitivity} this follows from the assumption that
  the Isomorphism Conjecture holds for every subgroup $H \subseteq
  G_i$ and in particular for any $H \in \psi_i^*\calc(G)$ for
  $\calc(G_i)|_H = \calc(H)$.
\end{proof}


\typeout{------------ Section 4: Fibered Isomorphism Conjectures and colimits
  ----------------}

\section{Fibered Isomorphism Conjectures and colimits}
\label{sec:Fibered_Isomorphism_Conjectures_and_colimits}

In this section we also deal with the Fibered version of the Isomorphism Conjectures.
(This is not directly needed for the purpose of this paper and the reader may skip this
section.)  This is a stronger version of the Farrell-Jones Conjecture. The Fibered
Farrell-Jones Conjecture does imply the Farrell-Jones Conjecture and has better
inheritance properties than the Farrell-Jones Conjecture.

We generalize (and shorten the proof of) the result of
Farrell-Linnell~\cite[Theorem~7.1]{Farrell-Linnell(2003a)} to a more general setting about
equivariant homology theories as developed in
Bartels-L\"uck~\cite{Bartels-Lueck(2004ind)}.

\begin{definition}[Fibered Isomorphism Conjecture for $\calh^?_*$]
\label{def:(Fibered)_Isomorphism_Conjectures_for_calh?_ast}
Fix a group $\Gamma$ and an equivariant homology theory $\calh^?_*$ with values in
$\Lambda$-modules over $\Gamma$.  A group $(G,\xi)$ over $\Gamma$ together with a family
of subgroups $\calf$ of $G$ satisfies the \emph{Fibered Isomorphism Conjecture (for
  $\calh^?_*$)} if for each group homomorphism $\phi \colon K \to G$ the group $(K,\xi
\circ \phi)$ over $\Gamma$ satisfies the Isomorphism Conjecture with respect to the family
$\phi^*\calf$.
\end{definition}

\begin{theorem}
\label{the:FIC_and_limits}
Let $(G,\xi)$ be a group over $\Gamma$. Let $\calf$ be a family of subgroups of $G$. Let
$\{G_i \mid i \in I\}$ be a directed system of groups with $G = \colim_{i \in I} G_i$ and
structure maps $\psi_i \colon G_i \to G$.  Suppose that $\calh^?_*$ is strongly continuous
and for every $i \in I$ the Fibered Isomorphism Conjecture holds for $(G_i,\xi
\circ\psi_i)$ and $\psi_i^*\calf$.

Then the Fibered Isomorphism Conjecture holds for $(G,\xi)$ and $\calf$.

\end{theorem}

\begin{proof} Let $\mu \colon K \to G$ be a group homomorphism. Consider the pullback of groups
  $$\begin{CD} K_i @ > \mu_i >> G_i
    \\
    @V \overline{\psi_i} VV @V \psi_i VV
    \\
    K @> \mu >> G
\end{CD}
$$
Explicitly $K_i = \{(k,g_i) \in K \times G_i \mid \mu(k) = \psi_i(g_i)\}$. Let
$\overline{\phi}_{i,j} \colon K_i \to K_j$ be the map induced by $\phi_{i,j} \colon G_i
\to G_j$, $\id_K$ and $\id_G$ and the pullback property. One easily checks by inspecting
the standard model for the colimit over a directed set that we obtain a directed system
$\overline{\phi}_{i,j} \colon K_i \to K_j$ of groups indexed by the directed set $I$ and
the system of maps $\overline{\psi}_i \colon K_i \to K$ yields an isomorphism $\colim_{i
  \in I} K_i \xrightarrow{\cong} K$. The following diagram commutes
$$\begin{CD} \colim_{i \in I}
  \calh_n^{K_i}\left(\overline{\psi_i}^*\mu^*\EGF{G}{\calf}\right) @>
  t^{K}_n(\mu^*\EGF{G}{\calf}) > \cong > \calh_n^{K}\left(\mu^*\EGF{G}{\calf}\right)
  \\
  @VVV @VVV
  \\
  \colim_{i \in I} \calh_n^{K_i}(\pt) @> t^{K}_n(\pt)> \cong> \calh_n^{K}(\pt)
\end{CD}
$$
where the vertical arrows are induced by the obvious projections onto $\pt$ and the
horizontal maps are the isomorphisms from Lemma~\ref{lem:directed_limits_and_calh?_ast}.
Notice that $\overline{\psi_i}^*\mu^*\EGF{G}{\calf}$ is a model for
$\EGF{K_i}{\overline{\psi_i}^*\mu^*\calf} = \EGF{K_i}{\mu_i^*\psi_i^*\calf}$. Hence each
map $\calh_n^{K_i}\left(\overline{\psi_i}^*\mu^*\EGF{G}{\calf}\right) \to
\calh_n^{K_i}(\pt)$ is bijective since $(G_i,\xi \circ \psi_i)$ satisfies the Fibered
Isomorphism Conjecture for $\psi_i^*\calf$ and hence $(K_i,\xi \circ \psi_i \circ \mu_i)$
satisfies the Isomorphism Conjecture for $\mu_i^*\psi_i^*\calf$. This implies that the
left vertical arrow is bijective. Hence the right vertical arrow is an isomorphism. Since
$\mu^*\EGF{G}{\calf}$ is a model for $\EGF{K}{\mu^*\calf}$, this means that $(K,\xi \circ
\mu)$ satisfies the Isomorphism Conjecture for $\mu^*\calf$. Since $\mu \colon K \to G$ is
any group homomorphism, $(G,\xi)$ satisfies the Fibered Isomorphism Conjecture for
$\calf$.
\end{proof}

The proof of the following results are analogous to the one
in~\cite[Lemma~1.6]{Bartels-Lueck(2004ind)} and~\cite[Lemma~1.2]{Bartels-Lueck(2006)},
where only the case $\Gamma = \{1\}$ is treated.

\begin{lemma} \label{lem:enlarging_the_family_for_the_Fibered_Isomorphism_Conjecture}
  Let $(G,\xi)$ be a group over $\Gamma$ and let $\calf \subset \calg$ be families of
  subgroups of $G$.  Suppose that $(G,\xi)$ satisfies the Fibered Isomorphism Conjecture
  for the family $\calf$.

  Then $(G,\xi)$ satisfies the Fibered Isomorphism Conjecture for the family $\calg$.
\end{lemma}

\begin{lemma} \label{lem:basic_inheritance_property_of_fibered_conjecture}
  Let $(G,\xi)$ be a group over $\Gamma$. Let $\phi \colon K \to G$ be a group
  homomorphism and let $\calf$ be a family of subgroups of $G$.  If $(G,\xi)$ satisfies
  the Fibered Isomorphism Conjecture for the family $\calf$, then $(K,\xi \circ \phi)$
  satisfies the Fibered Isomorphism Conjecture for the family $\phi^*\calf$.
\end{lemma}

For the remainder of this section fix a class of groups $\calc$ closed under isomorphisms,
taking subgroups and taking quotients, e.g., the families $\Fin$ or $\VCyc$.

\begin{lemma} \label{lem:Fibered_Conjecture_passes_to_subgroups}
  Let $(G,\xi)$ be a group over $\Gamma$.  Suppose that the Fibered Isomorphism Conjecture
  holds for $(G,\xi)$ and $\calc(G)$. Let $H \subseteq G$ be a subgroup.

  Then the Fibered Isomorphism Conjecture holds for $(H,\xi|_H)$ and $\calc(H)$.
\end{lemma}

\begin{proof}
  This follows from Lemma~\ref{lem:basic_inheritance_property_of_fibered_conjecture}
  applied to the inclusion $H \to G$ since $\calc(H) = \calc(G)|_H$.
\end{proof}

\begin{theorem}
\label{the:fibered_isomorphism_conjecture_is_stable_under_colim}
Let $(G,\xi)$ be a group over $\Gamma$.
\begin{enumerate}
\item
\label{the:fibered_isomorphism_conjecture_is_stable_under_colim:injective}
Let $G$ be the directed union $G = \bigcup_{i \in I} G_i$ of subgroups $G_i$. Suppose that
$\calh^?_*$ is continuous and that the Fibered Isomorphism Conjecture is true for
$(G_i,\xi|_{G_i})$ and $\calc(G_i)$ for all $i \in I$.

Then the Fibered Isomorphism Conjecture is true for $(G,\xi)$ and $\calc(G)$;

\item
\label{the:fibered_isomorphism_conjecture_is_stable_under_colim:general}
Let $\{G_i \mid i \in I\}$ be a directed system of groups with $G = \colim_{i \in I} G_i$
and structure maps $\psi_i \colon G_i \to G$.  Suppose that $\calh^?_*$ is strongly
continuous and that the Fibered Isomorphism Conjecture is true for $(G_i,\xi \circ
\psi_i)$ and $\calc(G_i))$ for all $i \in I$.

Then the Fibered Isomorphism Conjecture is true for $(G,\xi)$ and $\calc(G)$.
\end{enumerate}
\end{theorem}
\begin{proof}\ref{the:isomorphism_conjecture_is_stable_under_colim:injective}
  The proof is analogous to the one in~\cite[Proposition~3.4]{Bartels-Lueck(2006)}, where
  the case $\Gamma = \{1\}$ is considered.
  \\[1mm]\ref{the:isomorphism_conjecture_is_stable_under_colim:general} Because $\calc$ is
  closed under taking quotients we conclude $\calc(G_i) \subseteq \psi_i^*\calc(G)$.  Now
  the claim follows from Theorem~\ref{the:FIC_and_limits} and
  Lemma~\ref{lem:enlarging_the_family_for_the_Fibered_Isomorphism_Conjecture}.
\end{proof}

\begin{corollary} \label{cor:reducing_the_IC}

\begin{enumerate}

\item \label{cor:reducing_the_IC:unfibered}
      Suppose that $\calh^?_*$ is continuous. Then the (Fibered) Isomorphism
      Conjecture for $(G,\xi)$ and $\calc(G)$ is true for all groups $(G,\xi)$ 
      over $\Gamma$ if
      and only if it is true for all such groups where $G$ is a finitely generated group;

\item \label{cor:reducing_the_IC:fibered} Suppose that $\calh^?_*$ is strongly
      continuous. Then the Fibered Isomorphism Conjecture for $(G,\xi)$ and $\calc(G)$
      is true for all groups $(G,\xi)$ over $\Gamma$ if and only if it is true for all
      such groups where $G$ is finitely presented.

\end{enumerate}

\end{corollary}
\begin{proof}
  Let $(G,\xi)$ be a group over $\Gamma$ where $G$ is finitely generated. 
  Choose a finitely generated
  free group $F$ together with an epimorphism $\psi \colon F \to G$.
  Let $K$ be the kernel of $\psi$. Consider the directed system of
  finitely generated subgroups $\{K_i \mid i \in I\}$ of $K$. Let
  $\overline{K_i}$ be the smallest normal subgroup of $K$ containing
  $K_i$.  Explicitly $\overline{K_i}$ is given by elements which can
  be written as finite products of elements of the shape $fk_if^{-1}$
  for $f \in F$ and $k \in K_i$. We obtain a directed system of groups
  $\{F/\overline{K_i} \mid i \in I\}$, where for $i \le j$ the
  structure map $\phi_{i,j} \colon F/\overline{K_i} \to
  F/\overline{K_j}$ is the canonical projection. If $\psi_i \colon
  F/\overline{K_i} \to F/K = G$ is the canonical projection, then the
  collection of maps $\{\psi_i \mid i \in I\}$ induces an
  isomorphism $\colim_{i \in I} F/\overline{K_i} \xrightarrow{\cong}
  G$. By construction for each $i \in I$ the group $F/\overline{K_i}$ is finitely
  presented and the Fibered Isomorphism Conjecture holds for 
  $(F/\overline{K_i},\xi \circ \psi_i)$ and $\calc(F/\overline{K_i})$ by
  assumption.
  Theorem~\ref{the:fibered_isomorphism_conjecture_is_stable_under_colim}
  ~\ref{the:fibered_isomorphism_conjecture_is_stable_under_colim:general}
  implies that the Fibered Farrell-Jones Conjecture for $(G,\xi)$ and $\calc(G)$ is true.
\end{proof}


\typeout{-------------------- Section 5: Some equivariant homology theories ------------}

\section{Some equivariant homology theories}
\label{sec:Some_equivariant_homology_theories}

In this section we will describe the relevant homology theories
over a group $\Gamma$ and show that they are (strongly)
continuous. (We have defined the notion of an equivariant homology
theory over a group in
Definition~\ref{def:Equivariant_homology_theory_over_a_group_Gamma}.)


\subsection{Desired equivariant homology theories}
\label{subsec:desired_equivariant_homology_theories}
We will need the following

\begin{theorem}[Construction of equivariant homology theories]
\label{the:Construction_of_equivariant_homology_theories}
Suppose that we are given a group $\Gamma$ and a ring $R$ (with involution) or a
$C^*$-algebra $A$ respectively on which $\Gamma$ acts by structure preserving automorphisms.
Then:
\begin{enumerate}

\item \label{the:Construction_of_equivariant_homology_theories:list} Associated to
  these data there are equivariant homology theories with values in $\IZ$-modules
  over the group $\Gamma$
\begin{eqnarray*}
&& H_*^?(-;\bfK_R)
\\
&& H_*^?(-;\bfKH_R)
\\
&& H_*^?(-;\bfL_R^{\langle -\infty\rangle}),
\\
&& H_*^?(-;\bfK^{\topo}_{A,l^1}),
\\
&& H_*^?(-;\bfK^{\topo}_{A,r}),
\\
&& H_*^?(-;\bfK^{\topo}_{A,m}),
\end{eqnarray*}
where in the case $H_*^?(-;\bfK^{\topo}_{A,r})$ we will have to impose the restriction to the
induction structure that a homomorphisms $\alpha \colon (H,\xi) \to (G,\mu)$ over
$\Gamma$ induces a transformation $\ind_{\alpha}\colon \calh_n^H(X,X_0) \to
\calh_n^G(\alpha_*(X,X_0))$ only if the kernel of the underlying group 
homomorphism $\alpha \colon H \to G$
acts with amenable isotropy on $X\setminus X_0$;

\item \label{the:Construction_of_equivariant_homology_theories:coefficients}

  If $(G,\mu)$ is a group over $\Gamma$ and $H \subseteq G$ is a subgroup, then there are
  for every $n \in \IZ$ identifications
  $$\begin{array}{lclcl} H_n^H(\pt;\bfK_R) & \cong & H_n^G(G/H;\bfK_R) &
    \cong & K_n(R\rtimes H);
    \\
    H_n^H(\pt;\bfKH_R) & \cong & H_n^G(G/H;\bfKH_R) &
    \cong & \KH_n(R\rtimes H);
    \\
    H_n^H(\pt;\bfL_R^{\langle -\infty\rangle}) & \cong & H_n^G(G/H;\bfL_R^{\langle
      -\infty\rangle}) & \cong & L_n^{\langle-\infty\rangle}(R\rtimes H);
    \\
    H_n^H(\pt;\bfK^{\topo}_{A,l^1}) & \cong & H_n^G(G/H;\bfK^{\topo}_{A,l^1}) & \cong &
    K_n(A\rtimes_{l^1} H);
    \\
    H_n^H(\pt;\bfK^{\topo}_{A,r}) & \cong & H_n^G(G/H;\bfK^{\topo}_{A,r}) & \cong &
    K_n(A\rtimes_r H);
    \\
    H_n^H(\pt;\bfK^{\topo}_{A,m}) & \cong & H_n^G(G/H;\bfK^{\topo}_{A,r}) & \cong &
    K_n(A\rtimes_m H).
    \end{array}$$
    Here $H$ and $G$ act on $R$ and $A$ respectively via the given $\Gamma$-action,
    $\mu \colon G \to \Gamma$ and the inclusion $H \subseteq G$, $K_n(R\rtimes H)$ is the
    algebraic $K$-theory of the twisted group ring $R\rtimes H$, $\KH_n(R\rtimes H)$ is the
    homotopy $K$-theory of the twisted group ring $R\rtimes H$, $L^{\langle - \infty
    \rangle}_n(R\rtimes H)$ is the algebraic $L$-theory with decoration $\langle -\infty\rangle$
    of the twisted group ring with
    involution $R \rtimes H$, $K_n(A\rtimes_{l^1} H)$ is the topological $K$-theory of
    the crossed product Banach algebra $A\rtimes_{l^1} H$, $K_n(A\rtimes_r H)$ is the
    topological $K$-theory of the reduced crossed product $C^*$-algebra $A\rtimes_r H$,
    and $K_n(A\rtimes_m H)$ is the topological $K$-theory of the maximal crossed product
    $C^*$-algebra $A\rtimes_m H$;

\item \label{the:Construction_of_equivariant_homology_theories:Naturality_in_Gamma} Let
  $\zeta\colon \Gamma_0 \to \Gamma_1$ be a group homomorphism. Let $R$ be a ring (with
  involution) and $A$ be a $C^*$-algebra on which $\Gamma_1$ acts by structure preserving
  automorphisms. Let $(G,\mu)$ be a group over $\Gamma_0$. Then in all cases the evaluation
  at $(G,\mu)$ of the equivariant homology theory over $\Gamma_0$ associated to $\zeta^*R$ or
  $\zeta^*A$ respectively agrees with the evaluation at $(G,\zeta \circ \mu)$ of
  the equivariant homology theory over $\Gamma_1$ associated to $R$ or $A$ respectively.

\item \label{the:Construction_of_equivariant_homology_theories:Naturality_in_R_A} Suppose
  the group $\Gamma$ acts on the rings (with involution) $R$ and $S$ or on the
  $C^*$-algebras $A$ and $B$ respectively by structure preserving automorphisms. Let $\xi \colon R \to
  S$ or $\xi \colon A \to B$ be a $\Gamma$-equivariant homomorphism of rings (with involution)
  or $C^*$-algebras respectively. Then $\xi$ induces natural transformations of
  homology theories over $\Gamma$
$$\begin{array} {llcl}
\xi^?_* \colon &H_*^?(-;\bfK_R)
& \to &
H_*^?(-;\bfK_S);
\\
\xi^?_* \colon & H_*^?(-;\bfKH_R)
& \to &
H_*^?(-;\bfKH_S);
\\
\xi^?_* \colon & H_*^?(-;\bfL_R^{\langle -\infty\rangle})
& \to &
H_*^?(-;\bfL_S^{\langle -\infty\rangle});
\\
\xi^?_* \colon & H_*^?(-;\bfK^{\topo}_{A,l^1})
& \to &
H_*^?(-;\bfK^{\topo}_{B,l^1});
\\
\xi^?_* \colon & H_*^?(-;\bfK^{\topo}_{A,r})
& \to &
H_*^?(-;\bfK^{\topo}_{B,r});
\\
\xi^?_* \colon & H_*^?(-;\bfK^{\topo}_{A,m})
& \to &
H_*^?(-;\bfK^{\topo}_{B,m}).
\end{array}$$
They are compatible with the identifications appearing in
assertion~\ref{the:Construction_of_equivariant_homology_theories:coefficients};
\item
\label{the:Construction_of_equivariant_homology_theories_Naturality_in_K}
Let $\Gamma$ act  on the $C^*$-algebra $A$ by structure preserving
automorphisms. We can consider $A$ also as a ring with structure
preserving $G$-action. Then there are natural transformations of
equivariant  homology theories with values in $\IZ$-modules over
$\Gamma$
\begin{multline*} \hspace{20mm}
H_*^?(-;\bfK_A) \to H_*^?(-;\bfKH_A) \to H_*^?(-;\bfK^{\topo}_{A,l^1})
\\
\to H_*^?(-;\bfK^{\topo}_{A,m}) \to H_*^?(-;\bfK^{\topo}_{A,r}).
\end{multline*}
They are compatible with the identifications appearing in
assertion~\ref{the:Construction_of_equivariant_homology_theories:coefficients}.
\end{enumerate}
\end{theorem}


\subsection{(Strong) Continuity}
\label{subsec:(Strong)_Continuity}

Next we want to show

\begin{lemma} \label{lem:homology_theories_are_(strongly)_continuous}
  Suppose that we are given a group $\Gamma$ and a ring $R$ (with
  involution) or a $C^*$-algebra $A$ respectively on which $G$ acts by
  structure preserving automorphisms.

  Then the homology theories with values in $\IZ$-modules over
  $\Gamma$
  $$H_*^?(-;\bfK_R), H_*^?(-;\bfKH_R),
  H_*^?(-;\bfL_R^{\langle -\infty\rangle}),
  H_*^?(-;\bfK^{\topo}_{A,l^1}), \text{ and }
  H_*^?(-;\bfK^{\topo}_{A,m})$$
  (see
  Theorem~\ref{the:Construction_of_equivariant_homology_theories}) are
  strongly continuous in the sense of
  Definition~\ref{def:strongly_continuous_equivariant_homology_theory},
  whereas $$H_*^?(-;\bfK^{\topo}_{A,r})$$
  is only continuous.
\end{lemma}
\begin{proof} We begin with $H_*^?(-;\bfK_R)$ and $H_*^?(-;\bfKH_R)$.
  We have to show for every directed systems of groups $\{G_i \mid i
  \in I\}$ with $G = \colim_{i \in I} G_i$ together with a map $\mu
  \colon G \to \Gamma$ that the canonical maps
  \begin{eqnarray*}
  \colim_{i \in I} K_n(R\rtimes G_i)  & \to &  K_n(R\rtimes G);
  \\
  \colim_{i \in I} \KH_n(R\rtimes G_i)  & \to &  \KH_n(R\rtimes G),
  \end{eqnarray*}
  are bijective for all $n \in \IZ$.  Obviously $R\rtimes G$ is the
  colimit of rings $\colim_{i \in I} R\rtimes G_i$. Now the claim
  follows for $K_n(R\rtimes G)$ for $n \ge 0$ from \cite[(12) on
  page~20]{Quillen(1973)}.

  Using the Bass-Heller-Swan decomposition
  one gets the results for $K_n(R\rtimes G)$ for all $n \in \IZ$ and that
  the map
  \begin{eqnarray*}
  \colim_{i \in I} N^pK_n(R\rtimes G_i)  & \to &  N^pK_n(R\rtimes G)
  \end{eqnarray*}
  is bijective for all $n \in \IZ$ and all $p \in \IZ, p \ge 1 $ for
  the Nil-groups $N^pK_n(RG)$ defined by Bass~\cite[XII]{Bass(1968)}.
  Now the claim for homotopy $K$-theory follows from the spectral sequence due to
  Weibel~\cite[Theorem~1.3]{Weibel(1989)}.

  Next we treat $H_*^?(-;\bfL_R^{\langle -\infty\rangle})$.  We have
  to show for every directed systems of groups $\{G_i \mid i \in I\}$
  with $G = \colim_{i \in I} G_i$ together with a map $\mu \colon G
  \to \Gamma$ that the canonical map
  $$\colim_{i \in I} L_n^{\langle - \infty \rangle}(R\rtimes G_i)
  \to  L_n^{\langle - \infty \rangle}(R\times G)$$
  is bijective for all $n \in \IZ$.  Recall from \cite[Definition~17.1
  and Definition~17.7]{Ranicki(1992a)} that
  \begin{eqnarray*}
  L_n^{\langle-\infty\rangle}(R\rtimes G) & = & \colim_{m \to \infty}
  L_n^{\langle-m\rangle}(R\rtimes G);
  \\
  L_n^{\langle-m\rangle}(R\rtimes G) & = & \coker\left(L_{n+1}^{\langle-m+1\rangle}(R\rtimes G) \to
  L_{n+1}^{\langle-m+1\rangle}(R\rtimes G[\IZ])\right)  \; \text{ for } m \ge 0.
  \end{eqnarray*}
  Since $L_n^{\langle 1 \rangle}(R\rtimes G)$ is $L_n^h(R\rtimes G)$,
  it suffices to show that
  \begin{eqnarray}
  \omega_n \colon \colim_{i \in I} L_n^{h}(R\rtimes G_i) & \to & L_n^{h}(R\rtimes G) \label{omega_colim_L_to_L}
  \end{eqnarray}
  is bijective for all $n \in \IZ$. We give the proof of
  surjectivity for $n = 0$ only, the proofs of injectivity for $n = 0$ and
  of bijectivity for the other values of $n$ are similar.

  The ring $R\rtimes G$ is the colimit of rings $\colim_{i \in I}
  R\rtimes G_i$. Let $\psi_i \colon R\rtimes G_i \to R\rtimes G$ and
  $\phi_{i,j} \colon R\rtimes G_i \to R\rtimes G_j$ for $i,j \in I, i
  \le j$ be the structure maps.  One can define $R\rtimes G$ as the
  quotient of $\coprod_{i \in I} R\rtimes G_i/\sim$, where $x \in
  R\rtimes G_i$ and $y \in R\rtimes G_j$ satisfy $x \sim y$ if and
  only if $\phi_{i,k}(x) = \phi_{j,k}(y)$ holds for some $k \in I$
  with $i,j \le k$. The addition and multiplication is given by adding
  and multiplying representatives belonging to the same $R\rtimes
  G_i$.  Let $M(m,n;R\rtimes G)$ be the set of $(m,n)$-matrices with
  entries in $R\rtimes G$. Given $A_i \in M(m,n;R\rtimes G_i)$, define
  $\phi_{i,j}(A_i) \in M(m,n;R\rtimes G_j)$ and $\psi_i(A_i) \in
  M(m,n;R\rtimes G)$ by applying $\phi_{i,j}$ and $\psi_i$ to each
  entry of the matrix $A_i$. We need the following key properties
  which follow directly from inspecting the model for the colimit
  above:
  \begin{enumerate}
  \item Given $A \in M(m,n;R\rtimes G)$, there exists $i \in I$ and
    $A_i \in M(m,n;R\rtimes G_i)$ with $\psi_{i}(A_i) = A$;

  \item Given $A_i \in M(m,n;R\rtimes G_i)$ and $A_j \in
    M(m,n;R\rtimes G_j)$ with $\psi_i(A_i) = \psi_j(A_j)$, there
    exists $k \in I$ with $i,j \le k$ and $\phi_{i,k}(A_i) =
    \phi_{j,k}(A_j)$.
\end{enumerate}

An element $[A]$ in $L_0^h(R\rtimes G)$ is represented by a quadratic
form on a finitely generated free $R\rtimes G$-module, i.e., a matrix
$A \in GL_n(R\rtimes G)$ for which there exists a matrix $B \in
M(n,n;R\rtimes G)$ with $A = B + B^*$, where $B^*$ is given by
transposing the matrix $B$ and applying the involution of $R$
elementwise.  Fix such a choice of a matrix $B$. Choose $i \in I$ and
$B_i \in M(n,n;R\rtimes G_i)$ with $\psi_i(B_i) = B$.  Then
$\psi_i(B_i + B_i^*) = A$ is invertible. Hence we can find $j \in I$
with $i \le j$ such that $A_j := \phi_{i,j}(B_i + B_j^*)$ is
invertible. Put $B_j = \phi_{i,j}(B_i)$. Then $A_j = B_j + B_j^*$ and
$\psi_j(A_j) = A$. Hence $A_j$ defines an element $[A_j] \in
L_n^h(R\rtimes G_j)$ which is mapped to $[A]$ under the homomorphism
$L_n^h(R\rtimes G_j) \to L_n(R\rtimes G)$ induced by $\psi_j$. Hence
the map $\omega_0$ of \eqref{omega_colim_L_to_L} is surjective.

Next we deal with $H_*^?(-;\bfK^{\topo}_{A,l^1})$.
We have to show for every directed systems of groups $\{G_i \mid i \in I\}$
with $G = \colim_{i \in I} G_i$ together with a map $\mu \colon G
\to \Gamma$ that the canonical map
$$\colim_{i \in I} K_n(A \rtimes_{l^1} G_i) \to K_n(A \rtimes_{l^1} G)$$
is bijective for all $n \in \IZ$. Since topological $K$-theory is
a continuous functor, it suffices to show that the colimit (or
sometimes also called inductive limit) of the system of Banach
algebras $\{A \rtimes_{l^1}G_i\mid i \in I\}$ in the category of
Banach algebras with norm decreasing homomorphisms is $A
\rtimes_{l^1} G$. So we have to show that for any Banach algebra
$B$ and any system of (norm deceasing) homomorphisms of Banach algebras $\alpha_i
\colon A \rtimes_{l^1} G_i \to B$ compatible with the structure
maps $A \rtimes_{l^1} \phi_{i,j} \colon A \rtimes_{l^1} G_i \to A
\rtimes_{l^1} G_j$ there exists precisely one homomorphism of
Banach algebras $\alpha \colon A \rtimes_{l^1} G \to B$ with the
property that its composition with the structure map $A
\rtimes_{l^1} \psi_i \colon A \rtimes_{l^1} G_i \to A
\rtimes_{l^1} G$ is $\alpha_i$ for $i \in I$.

It is easy to see that in the category of  $\IC$-algebras
the colimit of the system $\{A \rtimes G_i\mid i \in I\}$ is $A
\rtimes G$ with structure maps $A \rtimes \psi_i \colon A \rtimes
G_i \to A \rtimes G$. Hence the restrictions of the homomorphisms
$\alpha_i$ to the subalgebras $A \rtimes G_i$ yields a
homomorphism of central $\IC$-algebras $\alpha' \colon A \rtimes G
\to B$ uniquely determined by the property that the composition of
$\alpha'$ with the structure map $A \rtimes \psi_i \colon A
\rtimes G_i \to A \rtimes G$ is $\alpha_i|_{A \rtimes G_i}$ for $i
\in I$. If $\alpha$ exists, its restriction to the dense
subalgebra $A \rtimes G$ has to be $\alpha'$. Hence $\alpha$ is
unique if it exists. Of course we want to define $\alpha$ to be
the extension of $\alpha'$ to the completion $A \rtimes_{l^1} G$
of $A \rtimes G$ with respect to the $l^1$-norm. So it remains to
show that $\alpha' \colon A \rtimes G \to B$ is norm decreasing.
Consider an element $ u \in A \rtimes G$ which is given by a
finite formal sum $u = \sum_{g\in F} a_g \cdot g$, where $F\subset G$ is some finite subset 
of $G$ and  
$a_g \in A$ for $g\in F$.
We can choose an index $j \in I$ and a finite set $F'\subset G_j$ such that
 $\psi_j|_{F'}:F'\to F$ is one-to-one. For $g\in F$ let $g'\in F'$ 
denote the inverse image of $g$ under this map.
 Consider the element $v = \sum_{g'\in F'} a_g \cdot g'$ in $A \rtimes G_j$. 
By construction we have $A \rtimes \psi_j(v) =u$ and $||v|| = ||u|| = \sum_{i=1}^n ||a_i||$.  
We conclude
$$
||\alpha'(u)|| = ||\alpha' \circ (A \rtimes\psi_j)(v)|| = || \alpha_j(v)||
 \le ||v|| = ||u||.
$$
The proof for $H_*^?(-;\bfK^{\topo}_{A,m})$ follows similarly, using the fact that 
by definition of the norm on $A\rtimes_mG$ every $*$-homomorphism of 
$A\rtimes G$ into a $C^*$-algebra $B$ extends uiniquely to $A\rtimes_mG$.
The proof for the continuity of $H_*^?(-;\bfK^{\topo}_{A,r})$ follows from
\cite[Theorem 4.1]{Baum-Millington-Plymen(2003)}.
\end{proof}

Notice that we have proved all promised results of the introduction as soon as we have
completed the proof of Theorem~\ref{the:Construction_of_equivariant_homology_theories}
which we have used as a black box so far.

\typeout{--- Section 6: From spectra over groupoids to equivariant homology theories ---}

\section{From spectra over groupoids to equivariant homology theories}
\label{sec:From_spectra_over_groupoids_to_equivariant_homology_theories}

In this section we explain how one can construct equivariant homology theories
from spectra over groupoids.

A \emph{spectrum} $\bfE = \{(E(n),\sigma(n)) \mid n \in \IZ\}$ is a
sequence of pointed spaces $\{E(n) \mid n \in \IZ\}$ together with
pointed maps called \emph{structure maps} $\sigma(n)\colon E(n) \wedge
S^1 \longrightarrow E(n+1)$.  A \emph{(strong) map} of spectra
(sometimes also called function in the literature) $\bff \colon \bfE
\to \bfE^{\prime}$ is a sequence of maps $f(n) \colon E(n) \to
E^{\prime}(n)$ which are compatible with the structure maps
$\sigma(n)$, i.e., we have $f(n+1) \circ \sigma(n) =
\sigma^{\prime}(n) \circ \left(f(n) \wedge \id_{S^1}\right)$ for all
$n \in \IZ$. This should not be confused with the notion of a map of
spectra in the stable category (see \cite[III.2.]{Adams(1974)}).
Recall that the homotopy groups of a spectrum are defined by
$$\pi_i(\bfE) := \colim_{k \to \infty} \pi_{i+k}(E(k)),$$
where the
system $\pi_{i+k}(E(k))$ is given by the composition
$$\pi_{i+k}(E(k)) \xrightarrow{S} \pi_{i+k+1}(E(k)\wedge S^1)
\xrightarrow{\sigma (k)_*} \pi_{i+k+1}(E(k +1))$$
of the suspension
homomorphism and the homomorphism induced by the structure map. We
denote by $\Spectra$ the category of spectra.

A \emph{weak equivalence} of spectra is a map $\mathbf{f}\colon \bfE
\to \bfF$ of spectra inducing an isomorphism on all homotopy groups.

Given a small groupoid $\calg$, denote by $\Groupoids\downarrow\calg$
the category of small groupoids over $\calg$, i.e., an object is a
functor $F_0 \colon \calg_0 \to \calg$ with a small groupoid as source
and a morphism from $F_0 \colon \calg_0 \to \calg$ to $F_1 \colon
\calg_1 \to \calg$ is a functor $F \colon \calg_0 \to \calg_1$
satisfying $F_1 \circ F = F_0$. We will consider a group $\Gamma$ as a
groupoid with one object and $\Gamma$ as set of morphisms.  An
\emph{equivalence} $F \colon \calg_0 \to \calg_1$ of groupoids is a functor
of groupoids $F$ for which there exists a functor of groupoids $F'
\colon \calg_1 \to \calg_0$ such that $F' \circ F$ and $F \circ F'$
are naturally equivalent to the identity functor.  A functor $F \colon
\calg_0 \to \calg_1$ of small groupoids is an equivalence of groupoids if
and only if it induces a bijection between the isomorphism classes of
objects and for any object $x \in \calg_0$ the map $\aut_{\calg_0}(x)
\to \aut_{\calg_1}(F(x))$ induced by $F$ is an isomorphism of groups.

\begin{lemma}\label{lem:from_spectra_over_groupoids_equiv_homology_theories}
Let $\Gamma$ be a group. Consider a covariant functor
$$\bfE \colon \Groupoids \downarrow \Gamma \to \Spectra$$
which sends equivalences of groupoids to weak equivalences of spectra.

Then we can associate to it an equivariant homology theory $\calh^?_*(,-;\bfE)$ (with
values in $\IZ$-modules) over $\Gamma$ such that for every group $(G,\mu)$ over $\Gamma$
and subgroup $H \subseteq G$ we have a natural identification
$$\calh_n^H(\pt;\bfE) = \calh_n^G(G/H,\bfE) = \pi_n(\bfE(H)).$$

If $\bfT \colon \bfE \to \bfF$ is a natural transformation of such functors $\Groupoids
\downarrow \Gamma \to \Spectra$, then it induces a transformation of equivariant homology
theories over $\Gamma$
$$\calh^?_*(-;\bfT) \colon \calh^?_*(,-;\bfE) \to \calh^?_*(,-;\bfF)$$
such that for every
group $(G,\mu)$ over $\Gamma$ and subgroup $H \subseteq G$ the homomorphism
$\calh_n^H(\pt;\bfT) \colon \calh_n^H(\pt;\bfE) \to \calh_n^H(\pt;\bfF)$ agrees under the
identification above with $\pi_n(\bfT(H)) \colon \pi_n(\bfE(H)) \to \pi_n(\bfF(H))$.
\end{lemma}
\begin{proof}
  We begin with explaining how we can associate to a group $(G,\mu)$ over $\Gamma$ a
  $G$-homology theory $\calh^G_*(-;\bfE)$ with the property that for every subgroup $H
  \subseteq G$ we have an identification
  $$\calh_n^G(G/H,\bfE) = \pi_n(\bfE(H)).$$
  We just follow the construction
  in~\cite[Section~4]{Davis-Lueck(1998)}. Let $\Or(G)$ be the \emph{orbit category} of
  $G$, i.e., objects are homogenous spaces $G/H$ and morphisms are $G$-maps. Given a
  $G$-set $S$, the associated \emph{transport groupoid} $t^G(S)$ has $S$ as set of objects
  and the set of morphisms from $s_0 \in S$ to $s_1 \in S$ consists of the subset $\{g \in
  G \mid gs_1=s_2\}$ of $G$. Composition is given by the group multiplication.  A $G$-map
  of sets induces a functor between the associated transport groupoids in the obvious way.
  In particular the projection $G/H  \to G/G$ induces a functor of groupoids $\pr_S \colon
  t^G(S) \to t^G(G/G) = G$.  Thus $t^G(S)$ becomes an object in $\Groupoids \downarrow
  \Gamma$ by the composite $\mu \circ \pr_S$. We obtain a covariant functor $t^G \colon
  \Or(G) \to \Groupoids \downarrow \Gamma$. Its composition with the given functor $\bfE$
  yields a covariant functor
  $$\bfE^G := \bfE \circ t^G\colon \Or(G) \to \Spectra.$$
  Now define
  $$\calh^G_*(X,A;\bfE) := \calh^G_*(-;\bfE^G),$$
  where $\calh^G_*(-;\bfE^G)$ is the
  $G$-homology theory which is associated to $\bfE^G \colon \Or(G) \to \Spectra$ and
  defined in~\cite[Section~4 and~7]{Davis-Lueck(1998)}. Namely, if $X$ is a
  $G$-$CW$-complex, we can assign to it a contravariant functor $\map_G(G/?,X) \colon
  \Or(G) \to \Spaces$ sending $G/H$ to $\map_G(G/H,X) = X^H$ and put $\calh_n^G(X;\bfE^G)
  := \pi_n(\map_G(G/?,X)_+ \wedge_{\Or(G)} \bfE^G)$ for the spectrum $\map_G(G/?,X)_+
  \wedge_{\Or(G)} \bfE^G$ (which is denoted in~\cite{Davis-Lueck(1998)} by
  $\map_G(G/?,X)_+ \otimes_{\Or(G)} \bfE^G$).

  Next we have to explain the induction structure. Consider a group
  homomorphism $\alpha \colon (H,\xi) \to (G,\mu)$ of groups over
  $\Gamma$ and an $H$-$CW$-complex $X$. We have to construct a
  homomorphism
  $$\calh^H_n(X;\bfE) \to \calh^G_n(\alpha_* X;\bfE).$$
  This will be
  done by constructing a map of spectra
  $$\map_H(H/?,X)_+ \wedge_{\Or(H)} \bfE^H \to \map_G(G/?,\alpha_*
  X)_+ \wedge_{\Or(G)} \bfE^G.$$
  We follow the constructions in~\cite[Section~1]{Davis-Lueck(1998)}.  The homomorphism $\alpha$
  induces a covariant  functor  $\Or(\alpha)\colon  \Or(H) \to \Or(G)$
  by sending $H/L$ to $\alpha_*(H/L) = G/\alpha(L)$. Given a
  contravariant functor $Y \colon \Or(H) \to \Spaces$, we can assign
  to it its induction with $\Or(\alpha)$ which is a contravariant functor
  $\alpha_*Y   \colon \Or(G) \to \Spaces$. Given a contravariant functor $Z \colon
  \Or(G) \to \Spaces$, we can assign to it its restriction which is
  the contravariant functor $\alpha^*Z := Z \circ \Or(\alpha) \colon \Or(G)
  \to \Spaces$. Induction $\alpha_*$ and $\alpha^*$ form an adjoint
  pair. Given an $H$-$CW$-complex $X$, there is a natural identification
  $\alpha_*\left(\map_H(H/?,X)\right) = \map_G(G/?,\alpha_*X)$. 
  Using \cite[Lemma~1.9]{Davis-Lueck(1998)}
  we get for an $H$-$CW$-complex $X$ a natural map of spectra
  $$\map_H(H/?,X)_+ \wedge_{\Or(H)} \alpha^*\bfE^G \to
  \map_G(G/?,\alpha_* X)_+ \wedge_{\Or(G)} \bfE^G.$$
  Given an $H$-set $S$, we obtain a functor of groupoids
  $t^H(S) \to t^G(\alpha_*S)$ sending $s \in S$ to $(1,s) \in G \times_{\alpha} S$ and
  a morphism in $t^H(S)$ given by a group element $h$ to the one in $t^G(\alpha_*S)$
  given by $\alpha(h)$. This yields  a natural transformation of  covariant functors
  $\Or(H) \to \Groupoids \downarrow \Gamma$ from $t^H \to t^G \circ \Or(\alpha)$.
  Composing with the functor $\bfE$ gives a natural transformation of  covariant functors
  $\Or(H) \to \Spectra$ from $\bfE^H$ to $\alpha^* \bfE^G$. 
  It induces a map of spectra
  $$\map_H(H/?,X)_+ \wedge_{\Or(H)} \bfE^H \to \map_H(H/?,X)_+ \wedge_{\Or(H)}
  \alpha^*\bfE^G.$$
  Its composition with the maps of spectra constructed beforehand yields the desired map
  of spectra $\map_H(H/?,X)_+ \otimes_{\Or(H)} \bfE^H \to \map_G(G/?,\alpha_* X)_+
  \otimes_{\Or(G)} \bfE^G$.

  We omit the straightforward proof that the axioms of an induction structure are satisfied.
  This finishes the proof of Theorem~\ref{the:Construction_of_equivariant_homology_theories}.

  The statement about the natural transformation $\bfT \colon \bfE \to \bfF$ is obvious.
\end{proof}


\typeout{-------------- Section 7: Some K-theory spectra associated to groupoids------}

\section{Some $K$-theory spectra associated to groupoids}
\label{sec:Some_K-theory_spectra_associated_to_groupoids}

The last step in completing the proof of
Theorem~\ref{the:Construction_of_equivariant_homology_theories} is to
prove the following
Theorem~\ref{the:necessary_functors_from_groupoids_to_spectra} because
then we can apply it in combination with
Lemma~\ref{lem:from_spectra_over_groupoids_equiv_homology_theories}.
(Actually we only need the version of
Theorem~\ref{the:necessary_functors_from_groupoids_to_spectra}, where
$\calg$ is given by a group $\Gamma$.)  Let
$\Groupoids^{\finker}\downarrow\calg$ be the subcategory of $
\Groupoids\downarrow\calg$ which has the same objects and for which a
morphism from $F_0 \colon \calg_0 \to \calg$ to $F_1 \colon \calg_1
\to \calg$ given by a functor $F \colon \calg_0 \to \calg_1$
satisfying $F_1 \circ F = F_0$ has the property that for every object
$x \in \calg_0$ the group homomorphism $\aut_{\calg_0}(x) \to
\aut_{\calg_1}(F(x))$ induced by $F$ has a finite kernel.
Denote by $\Rings$, $\astRings$,  and $\CastAlg$ the categories
of rings, rings with involution and $C^*$-algebras.

\begin{theorem} \label{the:necessary_functors_from_groupoids_to_spectra}
Let $\calg$ be a fixed groupoid. Let $R \colon \calg \to \Rings$, $R
\colon \calg \to \astRings$,  or $A
\colon \calg \to \CastAlg$ respectively be a covariant functor.
Then there exist covariant functors
\begin{alignat*}{4}
& \bfK_R \colon && \Groupoids\downarrow\calg && \to && \quad \Spectra;
\\
& \bfKH_R \colon && \Groupoids\downarrow\calg && \to && \quad \Spectra;
\\
& \bfL_R^{\langle -\infty \rangle} \colon && \Groupoids\downarrow\calg && \to && \quad \Spectra;
\\
& \bfK_{A,l^1}^{\topo}   \colon && \Groupoids\downarrow\calg && \to && \quad \Spectra;
\\
& \bfK_{A,r}^{\topo}    \colon  && \Groupoids^{\finker}\downarrow\calg \;&& \to && \quad \Spectra;
\\
& \bfK_{A,m}^{\topo}  \colon   && \Groupoids\downarrow\calg && \to && \quad \Spectra,
\end{alignat*}
together with natural transformations
\begin{alignat*}{4}
& \bfI_1 \colon && \bfK && \to && \quad\bfKH;
\\
& \bfI_2 \colon && \bfKH && \to && \quad \bfK_{A,l^1}^{\topo};
\\
& \bfI_3 \colon && \bfK_{A,l^1}^{\topo} && \to && \quad \bfK_{A,m}^{\topo};
\\
& \bfI_4 \colon && \bfK_{A,m}^{\topo} && \to && \quad \bfK_{A,r}^{\topo},
\end{alignat*}
of functors from $\Groupoids\downarrow\calg$ or $\Groupoids^{\finker}\downarrow\calg$
respectively to $\Spectra$ such that the following holds:

\begin{enumerate}

\item \label{the:K-theory_spectra_of_groupoids:weak_equivalences}
  Let $F_i \colon \calg_i \to \calg$ be objects for $i = 0,1$ and
  $F \colon F_0 \to F_1$ be a morphism between them in
  $\Groupoids\downarrow\calg$ or $\Groupoids^{\finker}\downarrow\calg$
  respectively such that the underlying functor of groupoids $F \colon
  \calg_0 \to \calg_1$ is an equivalence of groupoids. Then the
  functors send $F$ to a  weak equivalences of spectra;

\item \label{the:K-theory_spectra_of_groupoids:one_object}
  Let $F_0 \colon \calg_0 \to \calg$ be an object in
  $\Groupoids\downarrow\calg$ or $\Groupoids^{\finker}\downarrow\calg$
  respectively such that the underlying groupoid $\calg_0$ has only
  one object $x$. Let $G = \mor_{\calg_0}(x,x)$ be its automorphisms
  group. We obtain a ring $R(y)$, a ring $R(y)$ with involution,
  or a $C^*$-algebra $B(y)$ with $G$-operation by structure preserving maps from
  the evaluation of the functor $R$ or $A$ respectively at $y = F(x)$. Then:
  $$\begin{array}{lcl}
  \pi_n(\bfK_R(F)) & = & K_n(R(y)\rtimes G);
  \\
  \pi_n(\bfKH_R(F)) & = & \KH_n(R(y)\rtimes G);
  \\
  \pi_n(\bfL^{\langle-\infty\rangle}_R(F)) & = & L_n^{\langle-\infty\rangle}(R(y)\rtimes G);
  \\
  \pi_n(\bfK^{\topo}_{A(y),l^1}(F)) & = & K_n(A(y)\rtimes_{l^1} G);
  \\
  \pi_n(\bfK^{\topo}_{A(y),r}(F)) & = & K_n(A(y)\rtimes_r G);
  \\
  \pi_n(\bfK^{\topo}_{A(y),m}(F)) & = & K_n(A(y)\rtimes_m G),
  \end{array}$$
  where $K_n(R(y)\rtimes G)$ is the algebraic $K$-theory of the
  twisted group ring $R(y) \rtimes G$, $\KH_n(R(y)\rtimes G)$ is the homotopy $K$-theory of the
  twisted group ring $R(y) \rtimes G$, $L^{\langle - \infty \rangle}_n(R(y)\rtimes G)$ 
  is the algebraic $L$-theory with decoration $\langle -\infty\rangle$ 
  of the  twisted group ring  with involution $R(y) \rtimes G$,
  $K_n(A(y)\rtimes_{l^1} G)$ is the topological $K$-theory of the
  crossed product Banach algebra $A(y)\rtimes_{l^1} G$,
  $K_n(A(y)\rtimes_r G)$ is the topological $K$-theory of the reduced
  crossed product $C^*$-algebra $A(y)\rtimes_r G$, and
  $K_n(A(y)\rtimes_m G)$ is the topological $K$-theory of the maximal
  crossed product $C^*$-algebra $A(y)\rtimes_m G$.

  The natural transformations
  $\bfI_1$, $\bfI_2$, $\bfI_3$ and $\bfI_4$ become under this identifications
  the obvious change of rings and theory homomorphisms

\item \label{the:K-theory_spectra_of_groupoids:naturality}
  These constructions are in the obvious sense natural in $R$ and $A$ respectively and in $\calg$.
\end{enumerate}
\end{theorem}

We defer the details of the proof of Theorem~\ref{the:necessary_functors_from_groupoids_to_spectra}
in \cite{Bartels-Joachim-Lueck(2007)}. 
Its proof requires some work but there are many special cases
which have already been taken care of. If we would not insist on groupoids but
only on groups as input, these are the standard algebraic $K$- and $L$-theory spectra
or topological $K$-theory spectra associated to group rings, group Banach algebras and group $C^*$-algebras.
The construction for the algebraic $K$- and $L$-theory and the topological $K$-theory
in the case, where $G$ acts trivially on a ring $R$ or a $C^*$-algebra are already
carried out or can easily be derived from~\cite{Bartels-Lueck(2006)},~\cite{Davis-Lueck(1998)},  
and~\cite{Joachim(2003a)} except for the case of a Banach algebra. The case of the $K$-theory
spectrum associated to an additive category with $G$-action has already been carried out in
\cite{Bartels-Reich(2007)}. The main work which remains to do is to treat the
Banach case and to construct the relevant natural transformation from
$\bfKH$ to $\bfK_{A,l^1}^{\topo}$


\def\cprime{$'$} \def\polhk#1{\setbox0=\hbox{#1}{\ooalign{\hidewidth
  \lower1.5ex\hbox{`}\hidewidth\crcr\unhbox0}}}



\end{document}